\documentclass{amsart}

\usepackage{amsthm}
\usepackage{enumerate}
\usepackage{amsmath}
\usepackage{amssymb}
\usepackage{tikz}
\usepackage{amsfonts}
\usepackage{epsfig}
\usepackage{color,graphics}
\usepackage{hyperref}%

\usepackage[linesnumbered,vlined,boxed,fillcomment]{algorithm2e}
\SetKw{And}{ and }
\SetKw{All}{ all }
\SetKw{Or}{ or }
\SetKw{Break}{break}
\SetKw{Continue}{continue}
\SetKw{KwFrom}{from }
\makeatletter
\let\oldnl\nl
\newcommand{\nonl}{\renewcommand{\nl}{\let\nl\oldnl}}
\makeatother

\hypersetup{
  colorlinks=true,
  breaklinks=true,
  urlcolor= blue,
  linkcolor= blue,
  bookmarksopen=true,
  pdftitle={Slope factorization},
  pdfauthor={Adrien Poteaux},
}

\newtheorem{lem}{Lemma}
\newtheorem{cor}{Corollary}
\newtheorem{prop}{Proposition}
\newtheorem{defi}{Definition}
\newtheorem{thm}{Theorem}
\newtheorem{rem}{Remark}
\newtheorem{ex}{Example}

\newcommand{\Oc}{\mathcal{O}}
\newcommand{\F}{\mathcal{F}}

\newcommand{\rp}{\mathbb{R}}
\newcommand{\fp}{\mathbb{F}}

\newcommand{\qp}{\mathbb{Q}}

\newcommand{\np}{\mathbb{N}}

\newcommand{\lp}{\mathbb{L}}
\newcommand{\zp}{\mathbb{Z}}
\newcommand{\kp}{\mathbb{K}}
\newcommand{\ap}{\mathbb{A}}

\newcommand{\Ot}{\tilde{\mathcal{O}}}
\newcommand{\vl}{v_{\lambda}}
\newcommand{\vlp}{v_{\lambda'}}

\newcommand{\dl}{d_{\lambda}}

\newcommand{\tl}{\tau_{\lambda}}
\newcommand{\apl}{\mathbb{A}_{\lambda}}
\newcommand{\bpl}{\mathbb{B}_{\lambda}}
\newcommand{\tapl}{\tilde{\mathbb{A}}_{\lambda}}
\newcommand{\tbpl}{\tilde{\mathbb{B}}_{\lambda}}
\newcommand{\ml}{m_{\lambda}}
\newcommand{\mlp}{m_{\lambda'}}

\newcommand{\mlG}{m_{\lambda_G}}

\newcommand{\mlP}{m_{\lambda_P}}
\newcommand{\al}{a_{\lambda}}
\newcommand{\alp}{a_{\lambda'}}

\newcommand{\bl}{b_{\lambda}}
\newcommand{\blp}{b_{\lambda'}}

\newcommand{\iny}{\mathrm{in}_y}
\newcommand{\lty}{\mathrm{lt}_y}

\newcommand{\vlG}{v_{\lambda_G}}

\DeclareMathOperator{\ord}{ord}

\DeclareMathOperator{\lc}{lc}

\DeclareMathOperator{\Aut}{Aut}

\DeclareMathOperator{\Card}{Card}
\DeclareMathOperator{\Vol}{Vol}

\DeclareMathOperator{\Conv}{Conv}

\DeclareMathOperator{\D}{D}

\DeclareMathOperator{\In}{in}

\newcommand{\inl}{\In_\lambda}

\title{Improvements of convex-dense factorization of bivariate polynomials}

\author{Martin WEIMANN}

\begin{document}

\maketitle

\begin{abstract}
We develop a new algorithm for factoring a bivariate polynomial $F\in \kp[x,y]$ which takes fully advantage of the geometry of the Newton polygon of $F$. Under some non degeneracy hypothesis, the complexity is $\Ot(Vr_0^{\omega-1} )$ where $V$ is the volume of the polygon and $r_0$ is its minimal lower lattice length. The integer $r_0$ reflects some combinatorial constraints imposed by the polygon, giving a reasonable and easy-to-compute upper bound for the number of non trivial indecomposable Minkovski summands. The proof is based on a new fast factorization algorithm in $\kp[[x]][y]$ with respect to a slope valuation, a result  which has its own interest. 
\end{abstract}

\section{Introduction}


Factoring a bivariate polynomial $F\in \kp[x,y]$ over a field $\kp$ is a fundamental task of Computer Algebra which
received a particular attention since the years 1970s. We refer the reader to \cite[Chapter III]{GaGe13} and \cite{CheGalSurvey,CheLec,L2,L1} for a detailed historical account and an extended bibliography on the subject. For a dense polynomial of bidegree $(d_x,d_y)$, the current complexity is $\Oc(d_x d_y^{\omega})$ plus one univariate factorization of degree $d_y$ \cite{L2, L1}. Here, $2 \le \omega \le 3$  is so that we can multiply $n\times n$ matrices over $\kp$ with $\Oc(n^\omega)$ operations in $\kp$. The current theoretical bound is $\omega \approx 2.38$ \cite{GaGe13}, although $\omega$ is in practice closer to $3$  in most software implementations. 

\smallskip
In this paper, we will rather focus on finer complexity indicators attached to the Newton polygon $N(F)$, convex hull of the set of exponents of $F$. The polynomial $F$ is assumed to be represented by the list of its coefficients associated to the lattice points of $N(F)$, including zero coefficients. Following \cite{BL}, we talk about \emph{convex-dense} representation.  Assuming $N(F)$ two-dimensional, the size of $F$ can also be measured as the euclidean volume $V$ of $N(F)$ by Pick's formula. 

\smallskip
Various convex-dense factorization algorithms have been proposed in the last two decades, see  e.g. \cite{Gao,BL,Wei5,Wei6} and references therein.  In \cite{BL}, the authors compute in softly linear time a map $\tau\in \Aut(\zp^2)$ so that the volume of $\tau(N(F))$ is comparable to the volume of its bounding rectangle. Applying a classical dense algorithm on the resulting polynomial $\tau(F)$, they get a complexity estimate $\Oc(V n^{\omega-1})$ where $n$ is the width of the bounding rectangle, thus recovering the usual cost if $F$ is a dense polynomial. 
However, this algorithm does not take advantage of the combinatorial constraints imposed by Ostrowski's theorem, namely:
$$
N(GH)=N(G)+N(H)
$$
where $+$ indicates Minkowski sum. 
Regarding this issue, we developed in \cite{Wei5,Wei6} some convex-dense algorithms  based on toric geometry which take fully advantage of Ostrowski's combinatorial constraints. Unfortunately, the algorithm only works in characteristic zero and the complexity is not optimal.

In this note, we intend to show that under some non degeneracy hypothesis, it is in fact possible to take into account both the volume and Ostrowski's constraints, and so in arbitrary characteristic. Our complexity improves \cite{BL}, the gain being particularly significant when $N(F)$ has few Minkovski summands. 

\subsubsection*{Complexity model.} We work with computation trees \cite[Section
4.4]{BuClSh97}. We use an algebraic RAM model, counting only the number of arithmetic operations in $\kp$. We classically denote $\Oc()$ and $\Ot()$ to respectively hide constant
and logarithmic factors in our complexity results ; see
e.g. \cite[Chapter 25, Section 7]{GaGe13}. We use fast multiplication of polynomials, so that two polynomials in $\kp[x]$ of degree at most $d$ can be multiplied in softly linear time $\Ot(d)$.

\subsection{Fast convex-dense factorization}
Let $P\subset \rp^2$ be a lattice polygon. Let $\Lambda(P)$ be the \emph{lower boundary} of $P$, union of edges whose inward normal vectors have strictly positive second coordinate.  
The  \emph{(lower) lattice length} of $P$ is
$$
r(P):=\Card(\Lambda(P)\cap \zp^2)-1.
$$ 
As $r(PQ)=r(P)+r(Q)$, this integer gives an easy-to-compute upper bound for the number of indecomposable Minkovski summands of $P$ which are not a vertical segment (computing all Minkovski sum decompositions is NP-complete \cite{GL}). 

\medskip

Let $F=\sum c_{ij} x^j y^i\in \kp [x^{\pm 1},y^{\pm 1}]$. The support of $F$ is the set of exponents $(i,j)\in \zp^2$ such that $c_{ij}\ne 0$. Take care that the exponents of $y$ are represented by the horizontal axis. The Newton polygon $N(F)$ of $F$ is the convex hull of its support and we denote for short $\Lambda(F)$ its lower boundary.  

\begin{defi}\label{def:degenerated} 
We say that $F$ is not degenerated if for all edge $E\subset \Lambda(F)$, the \emph{edge polynomial} $y^{-\ord_y(F_E)} F_E$ is separable in $y$, where $F_E:=\sum_{(i,j)\in E\cap \zp^2} c_{ij} x^j y^i$.
\end{defi}

\medskip
Note that $F_E\in \kp[x^{\pm 1}][y]$ is quasi-homogeneous, hence its factorization reduces to a univariate factorization of degree the lattice length of $E$. 

\medskip
\noindent

Let us denote for short $V=\Vol(N(F))$ and $r=r(N(F))$. Note that $r\le d_y$. Due to Ostrowski's theorem, $r$ is an upper bound for the numbers of irreducible factors of $F$ of positive $y$-degree. Our main result is:

\begin{thm}\label{thm:main}
There exists a deterministic algorithm which given $F\in \kp[x,y]$ non degenerated, computes the irreducible factorization of $F$ over $\kp$ with 
\begin{enumerate}
\item $\Ot(r V)+\Oc(r^{\omega-1}V)$ operations in $\kp$ if $p=0$ or $p\ge 4V$, or
\item $\Ot(k r^{\omega-1} V)$ operations in $\fp_p$ if $\kp=\fp_{p^k}$, 
\end{enumerate}
plus some univariate factorizations over $\kp$ whose degree sum is $r$.
\end{thm}

As in \cite{BL}, we recover the usual complexity estimate $\Oc(d_x d_y^{\omega})$ when $F$ is a dense polynomial. However, Theorem \ref{thm:main} may improve significantly \cite{BL} when $F$ is non degenerated, as illustrated by the following example.

\begin{ex}\label{ex:Example 1}
Let $F$ of bidegree $(2n,2n)$, with Newton polygon 
$$
N(F)=\Conv((0,2),(2n,0),(0,2n),(2n,2n)).
$$
\noindent
The lower lattice length is $r=2$, which is a very strong combinatorial constraint: there is a unique Minkovski sum decomposition whose summands have positive volume (Figure \ref{fig1} below).

\begin{figure}[h]
	\caption{}
	\label{fig1}
	\scalebox{0.17}{
	\begin{picture}(0,0)%
\includegraphics{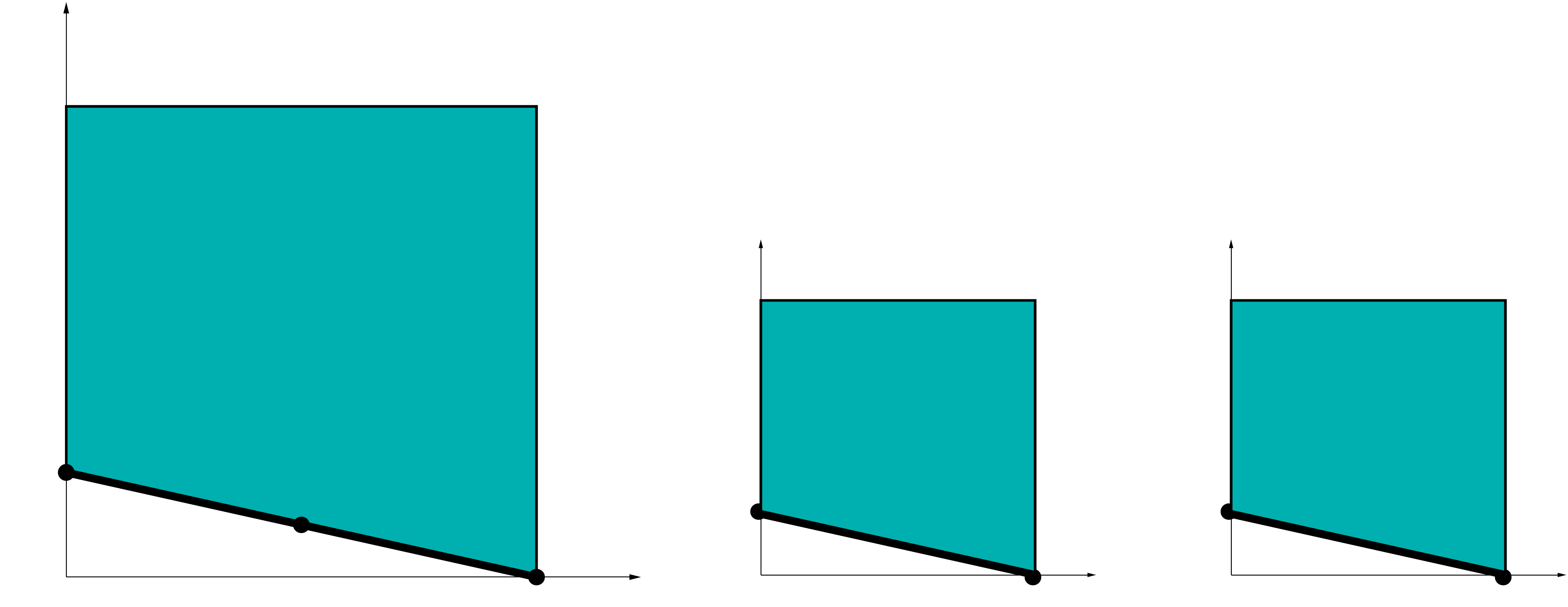}%
\end{picture}%
\setlength{\unitlength}{4144sp}%
\begingroup\makeatletter\ifx\SetFigFont\undefined%
\gdef\SetFigFont#1#2#3#4#5{%
  \reset@font\fontsize{#1}{#2pt}%
  \fontfamily{#3}\fontseries{#4}\fontshape{#5}%
  \selectfont}%
\fi\endgroup%
\begin{picture}(27000,10436)(-2489,-11825)
\put(6751,-12000){\scalebox{4}{$2n$}}
\put(-2200,-9511){\scalebox{4}{$2$}}
\put(-2274,-3211){\scalebox{4}{$2n$}}
\put(901,-10861){\scalebox{4}{$\Lambda$}}
\put(8551,-8611){\scalebox{4}{$=$}}
\put(16700,-8611){\scalebox{4}{$+$}}
\put(9901,-6586){\scalebox{4}{$n$}}
\put(15526,-11761){\scalebox{4}{$n$}}
\put(9901,-10186){\scalebox{4}{$1$}}
\put(17976,-10186){\scalebox{4}{$1$}}
\put(17976,-6586){\scalebox{4}{$n$}}
\put(23626,-11761){\scalebox{4}{$n$}}
\end{picture}%
}
\end{figure}

%

\noindent
As the bounding rectangle has size $\Oc(V)$, the convex-dense approach of \cite{BL} boils down to the dense algorithm \cite{L1}. We get the following complexity estimates:

\noindent

$\bullet$ Dense \cite{L1,L2} or convex-dense \cite{BL} algorithms: $\Oc(n^{\omega+1})$ operations in $\kp$ plus one univariate factorization of degree $2n$. 

\noindent

$\bullet$ Theorem \ref{thm:main} (assuming $F$ non degenerated): $\Ot(n^{2})$ operations in $\kp$ plus one univariate factorization of degree $2$.

\smallskip
\noindent
We get here a softly linear complexity. This is the most significant gain we can get, including the univariate factorization step. 
\end{ex}

A weakness of classical algorithms is to perform a shift $x\mapsto x+x_0$ to reduce to the case $F(0,y)$ separable, losing in such a way the combinatorial constraints offered by $N(F)$. Our approach avoids this shift.


\subsection{Even faster} We can play with affine automorphisms  $\tau\in \Aut(\zp^2)$ to minimize $r$ while keeping $V$ constant before applying Theorem \ref{thm:main}. This leads to the concept of \emph{minimal lattice length} of a lattice polygon $P$, defined as 
\begin{equation}
\label{eq:rmin}
r_0(P):=\min \{r(\tau(P))\,\, |\,\,  \tau\in \Aut(\zp^2)\}.
\end{equation}
This integer is easy to compute (Lemma \ref{lem:rmin}). Note that $r_{0}(N(F))$ can be reached by several $\tau$, which can lead to various lower boundaries with  lattice length $r_{0}$ (see Example \ref{ex:ex2} below). Let $\tau(F)$ be the image of $F$ when applying $\tau$ to its monomial exponents. 

\begin{defi}\label{def:minimally degenerated}  We say that $F$ is \emph{minimally non degenerated} if $\tau(F)$ is non degenerated for at least one transform $\tau$ reaching $r_{0}$.
\end{defi}

If $F$ is minimally non degenerated, we may apply Theorem \ref{thm:main} to $\tau(F)$, with same volume $V$ but with  smaller $r$. The factorization of $F$ is recovered for free from that of $\tau(F)$. We get:

\begin{cor}\label{cor:mmain}
Suppose that $F\in \kp[x,y]$ is minimally non degenerated with minimal lattice length $r_{0}$.  Then we can factorize $F$ with 
\begin{enumerate}
\item $\Ot(r_0 V)+\Oc(r_0^{\omega-1}V)$ operations in $\kp$ if $p=0$ or $p\ge 4V$, or
\item $\Oc(k r_{0}^{\omega-1} V)$ operations in $\fp_p$ if $\kp=\fp_{p^k}$, 
\end{enumerate}
plus some univariate factorizations over $\kp$ whose degree sum is $r_{0}$.
\end{cor}

Notice that similar transforms $F\mapsto \tau(F)$ are used in \cite{BL}, but the authors rather focus  on minimizing the size of the bounding rectangle of $N(F)$, while we focus on minimizing $r$. 
The following examples illustrate the differences between these two approaches. 

\begin{figure}[h]
	\caption{}
	\label{figComponent2}
	\bigskip
	\scalebox{0.3}{\begin{picture}(0,0)%
\includegraphics{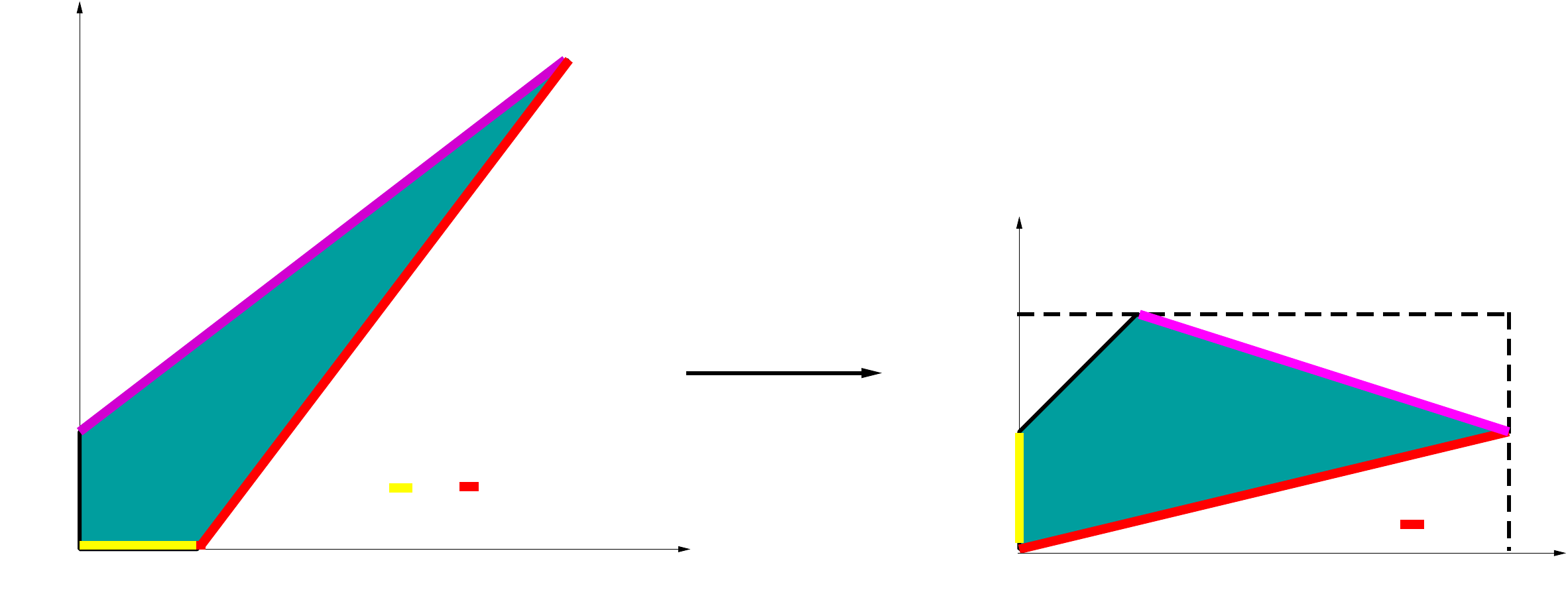}%
\end{picture}%
\setlength{\unitlength}{4144sp}%
\begingroup\makeatletter\ifx\SetFigFont\undefined%
\gdef\SetFigFont#1#2#3#4#5{%
  \reset@font\fontsize{#1}{#2pt}%
  \fontfamily{#3}\fontseries{#4}\fontshape{#5}%
  \selectfont}%
\fi\endgroup%
\begin{picture}(18012,6826)(211,-8450)
\put(2400,-8386){\scalebox{2.6}{$m$}}%
\put(6976,-2311){\scalebox{2.6}{$(n,n)$}}%
\put(526,-6586){\scalebox{2.6}{$m$}}%
\put(11326,-6586){\scalebox{2.6}{$m$}}
\put(12826,-4800){\scalebox{2.4}{$(m,2m)$}}
\put(17776,-6586){\scalebox{2.6}{$(n,m)$}}
\put(17551,-8386){\scalebox{2.6}{$n$}}
\put(11251,-5236){\scalebox{2.6}{$2m$}}
\put(9000,-5461){\scalebox{2.6}{$\tau$}}
\put(3676,-7291){\scalebox{2.6}{$\Lambda=$}}
\put(5101,-7291){\scalebox{2.6}{$+$}}
\put(14716,-7711){\scalebox{2.6}{$\Lambda_{\min}=$}}%
\end{picture}}
\end{figure}

\begin{ex}\label{ex:ex2} Let $0< m < n$ be two integers and suppose that 
$$
N(F)=\Conv((0,0),(m,0),(0,m),(n,n),
$$
as represented on the left side of Figure \ref{figComponent2}. 
The lower boundary $\Lambda(F)$ is the union of the yellow and red edges, with lattice length $r=m+\gcd(m,n)$. Applying the affine automorphism $\tau:(i,j)\mapsto (j,m-i+j)$, the resulting polygon $\tau(N(F))$ has red lower boundary, with minimal lattice length $r_{0}=gcd(m,n)$. The bounding rectangle of $\tau(N(F))$ has volume $2mn=V/2$, so \cite{BL} would apply a dense algorithm on $\tau(F)$. We get the following estimates:

\begin{itemize}
\item Dense algorithm \cite{L2,L1}: $\Oc(n^{\omega+1})$ operations and one univariate factorization of degree $n$.
\item Convex-dense algorithm  \cite{BL}: $\Oc(n m^{\omega})$ operations and one univariate factorization of degree $2m$.
\item Theorem \ref{thm:main} (assuming $F$ non degenerate): $\Oc(n m \gcd(n,m)^{\omega-1})$ operations and one univariate factorization of degree $\gcd(m,n)$.
\end{itemize}
Again, if $\gcd(m,n)\ll m$, our approach will be significantly faster than \cite{BL}, including the univariate factorization step. Notice that by symmetry, $r_{0}$ is reached also by the transform $\tau'$ which maps the purple edge as the lower convex hull. Hence, even if $F$ were "red-edge" degenerated, we would have a second chance that $F$ is not "purple-edge" degenerated, allowing then to apply Corollary \ref{cor:mmain}. 
\end{ex}

In the previous example, the image $\tau(F)$ reached simultaneously a minimal lower lattice length and a bounding rectangle of size $\Oc(V)$. The next example illustrates that this is not always the case.

\begin{figure}[h]
	\caption{}
	\label{figComponent3}
	\bigskip
\scalebox{0.2}{\begin{picture}(0,0)%
\includegraphics{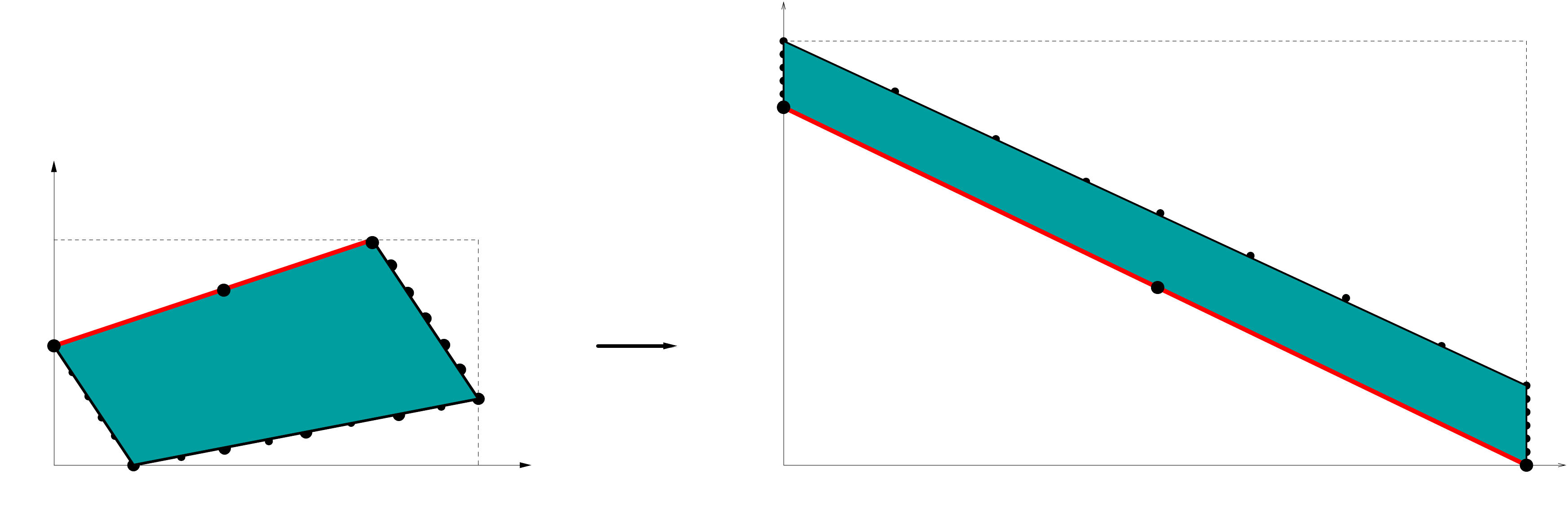}%
\end{picture}%
\setlength{\unitlength}{4144sp}%
\begingroup\makeatletter\ifx\SetFigFont\undefined%
\gdef\SetFigFont#1#2#3#4#5{%
  \reset@font\fontsize{#1}{#2pt}%
  \fontfamily{#3}\fontseries{#4}\fontshape{#5}%
  \selectfont}%
\fi\endgroup%
\begin{picture}(26577,8626)(12811,-14525)
\put(13051,-11761){\scalebox{3}{$2n$}}
\put(23476,-11311){\scalebox{4}{$\tau$}}
\put(12326,-9961){\scalebox{3}{$4n+4$}}
\put(16901,-9511){\scalebox{3}{$(kn-n-2,4n+4)$}}
\put(21151,-12661){\scalebox{3}{$(kn,n)$}}
\put(20701,-14236){\scalebox{3}{$kn$}}
\put(14851,-14236){\scalebox{3}{$n$}}
\put(25201,-6586){\scalebox{3}{$kn$}}
\put(24500,-7811){\scalebox{3}{$kn-n$}}
\put(38926,-12211){\scalebox{3}{$((2k-1)n,n+2)$}}
\put(38201,-14461){\scalebox{3}{$(2k-1)n$}}
\end{picture}%
}
\end{figure}

\begin{ex}\label{ex:ex3}
Suppose that $F$ has Newton polygon $N(F)$ as represented on the left side of figure \ref{figComponent3}, depending on parameters $k,n$. 
The bounding rectangle of $N(F)$ has volume $\Oc(kn^2)=\Oc(V)$, so \cite{BL} applies a dense algorithm on $F$. Any black edge has lattice length $n$ or $n+2$ while the red edge has lattice length  $r2$. We check that the affine automorphism $\tau(i,j)=(2i+j-2n,-i+kn)$ sends $N(F)$ to the right hand polygon, leading to $r_{0}=2$. We get the complexity estimates:
\begin{itemize}
\item  Dense \cite{L2,L1} or convex-dense algorithms  \cite{BL}: $\Oc(k n^{\omega+1})$ and one univariate factorization of degree $4n+4$.
\item  Theorem \ref{thm:main} (assuming $F$ minimally non degenerated): $\Ot(kn^2)$ operations and one univariate factorization of degree $2$. 
\end{itemize}

\noindent
Again, we get a softly linear complexity. This example illustrates the fact that minimizing the lower lattice length may increase significantly the volume of the bounding rectangle ($k^2 n^2 \gg V$). 
\end{ex}

Classical fast factorization algorithms are based on a "lifting and recombination" scheme: factorize $F$ in $\kp[[x]][y]$ with $x$-adic precision $\Oc(d_x)$ and recombine the analytic factors into global factors. Example \ref{ex:ex3} shows that we can not apply this strategy to our target polynomial $\tau(F)$: the analytic factorization with precision $d_x=kn$ would have size $\Oc(k^2 n^2)$ which does not fit in our aimed bound. To remediate this, we will rather factorize $\tau(F)$ in $\kp[[x]][y]$ with respect to another suitable valuation depending on the polygon. This is the second main result of our paper, that we explain now.

\subsection{Fast valuated analytic factorization.} Let $\lambda\in \qp$ and let $\vl$ stands for the valuation
\begin{equation}\label{eq:vlambda_intro}
v_\lambda:\kp((x))[y]\to \qp,\qquad v_\lambda(\sum c_{ij} x^j y^i):=\min(j+i\lambda, \,c_{ij}\ne 0),
\end{equation}
with convention $\vl(0)=\infty$. 
%
If $F \in \kp((x))[y]$, the lower convex hull $\Lambda=\Lambda(F)$ is well defined, and Definition \ref{def:degenerated} still makes sense in this larger ring. We denote 
\begin{equation}\label{eq:mlambda_intro}
\ml(F)=\max_{(i,j)\in \Lambda}(j+i\lambda) -\vl(F).
\end{equation}
Note that $\ml(F)\ge 0$, with equality if and only if $\Lambda(F)$ is straight of slope $-\lambda$. 
We measure the quality of the $\vl$-approximation of $F$ by a polynomial $G$ by the relative quantity $\vl(F-G)-\vl(F)$. We prove:

\begin{thm}\label{thm:lambda-facto}
Let $F\in \kp((x))[y]$ monic of degree $d$. Suppose that $F$ is non degenerate, with monic irreducible factors $F_1^*,\ldots,F_s^*$. Given $\sigma\ge \ml(F)$, we can compute $F_1,\ldots,F_s$ monic such that 
$$
\vl(F-F_1\cdots F_s)-\vl(F)> \sigma
$$  
with $\Ot(d\sigma)$ operations in $\kp$ plus some univariate factorizations over $\kp$ whose degree sum is at most $d$. Moreover, each factor is approximated with a relative precision
$$
\vl(F_i-F_i^*)- \vl(F_i^*) > \sigma-\ml(F).
$$
for all $i=1,\ldots,s$.
\end{thm}

Up to our knowledge, this result is new. It improves \cite{PoWe19} and \cite{PoWe22}, which focus on the Gauss valuation  $v_0$ and reach a quasi-optimal complexity only for $\sigma \ge d m_0(F)$ and characteristic of $\kp$ zero or high enough. It turns out that we need to get rid of all these restrictions for our purpose. The proof of Theorem \ref{thm:lambda-facto} is based on two main points:

$\bullet$ Fast arithmetic of sparse polynomials, leading to a softly linear $\vl$-adic Hensel lifting (Proposition \ref{prop:lambda-hensel}).

$\bullet$ A divide and conquer strategy based on a suitable choice of the various slopes $\lambda'$ which will be used at each recursive call of Hensel lifting.

\subsection{Main lines of the proof of Theorem \ref{thm:main}}

Except the choice of the valuation, the  strategy for the proof of Theorem \ref{thm:main} mainly follows \cite{L1,We16}:

\medskip

$\bullet$ We choose a suitable $\lambda\in \qp$ and we compute the factorization of $F$ in $\kp[[x]][y]$ with $\vl$-adic precision $\sigma\in \Oc(V/d_y)$, for a cost $\Ot(V)$ by Theorem \ref{thm:lambda-facto}. 

$\bullet$ Adapt the logarithmic derivative method of  \cite{L1, We16} to reduce to linear algebra the problem of recombinations of the truncated analytic factors into factors in $\kp[x,y]$.  \emph{A good choice of $\lambda$ is a key point} to ensure that the $\vl$-adic precision $\Oc(V/d_y)$ is sufficient to solve recombinations. 

$\bullet$ We are reduced to solve a linear system of at most $r$ unknowns and $\Oc(V)$ equations, which fits in the aimed bound. We build the underlying recombination matrix using a fast $\vl$-adic euclidean division by non monic polynomials (Proposition \ref{prop:final_division}). 

\begin{rem}
If $F$ is degenerated, we may probably compute nevertheless in softly linear time a $\vl$-adic factorization of $F\in \kp[[x]][y]$ using  recent algorithms  \cite{PoWe19,PoWe22} combined with Theorem \ref{thm:lambda-facto}. The number of factors to recombine is less than $r_{0}$, and possibly much smaller. Unfortunately, we might need a higher precision for solving recombinations, in which case the cost does not fit in the aimed bound. We refer the reader to \cite{We16} for mode details of such an approach in the $x$-adic case. 
\end{rem}

\begin{rem}
Let us mention too \cite{CaRoVa16}, where the authors develop a Hensel lifting with respect to a \textit{Newton precision}, given by a convex piecewise affine function. It might be interesting to look if such an approach could be useful for our purpose, as it allows to take care of the shape of $\Lambda(F)$.
\end{rem}

\subsection{Organisation of the paper} 
Section \ref{sec:analytic_facto} is dedicated to the proof of Theorem \ref{thm:lambda-facto}.
In section \ref{sec:facto}, we adapt the lifting and recombination scheme of \cite{L1,L2} in the $\vl$-adic context, leading to the proof of Theorem \ref{thm:main} and Corollary \ref{cor:mmain}.

\section{Fast $\vl$-adic factorization}\label{sec:analytic_facto}

In what follows, we fix $\lambda=m/q\in \qp$ with $q\ge 1$ and $q,m$ coprime and we consider the valuation $\vl$ as defined in \eqref{eq:vlambda_intro}.

\subsection{The ring $\apl$ and its fast arithmetic}\label{ssec:ring_apl}

Consider the classical Newton-Puiseux transformation
\begin{equation}\label{eq:tau_lambda}
\tau_\lambda:\kp((x))[y] \to \kp((x))[y], \qquad F(x,y)\mapsto \hat{F}(x,y)=F(x^q,x^m y).
\end{equation}
This map is an injective $\kp$-algebra endomorphism. Thus, its image 
$$
\apl:=\kp((x^q))[x^m y]\subset \kp((x))[y]
$$ 
is a subring isomorphic to  $\kp((x))[y]$. We denote
$$
\apl^+=\apl\cap \kp[[x]][y]\qquad {\rm and}\qquad \bpl=\apl\cap \kp[x,y]
$$
Both sets are subrings of $\apl$. Note that the map $\tau_\lambda$ preserves the size of the support of a polynomial. 

The valuation $\vl$ is related to the Gauss valuation $v_0$ by
\begin{equation}\label{eq:v0vl}
v_0(\tau_\lambda(F))=q \vl(F)
\end{equation}
Unfortunately, computing the $v_0$-adic factorization of $\tau_\lambda(F)$ which induces the $\vl$-adic factorization of $F$ with the recent softly linear algorithms \cite{PoWe20} does not fit in the aimed bound due to the presence of the extra factor $q$ in \eqref{eq:v0vl}. To remediate this problem, we need to take advantage of the fact that $\tau_\lambda(F)$ is sparse, which is reflected in more details by the following lemma:

\begin{lem}\label{lem:the_ring_apl}
Let $F\in \kp((x))[y]$. Then $F\in \apl$ if and only if
$$
F=\sum_{k} f_k(y^q) y^{\alpha_\lambda(k)}  x^k,\quad f_k\in \kp[y]
$$
where $0\le \alpha_\lambda(k)< q$ is defined by $
\alpha_\lambda(k)\equiv k\,m^{-1} \mod q.$
Equivalently, we have:
$$
\apl=\bigoplus_{k=0}^{q-1} x^k y^{\alpha_\lambda(k)} \,\kp((x^q))[y^q].
$$
In particular, $\apl\cap \kp((x))=\kp((x^q))$ and $ \apl\cap\kp[y]=\kp[y^q].$
\end{lem}

\begin{proof}
By \eqref{eq:tau_lambda}, we have $F\in \apl$ if and only if $F=\sum_{i,j} c_{ij} x^{mi+qj} y^i$ for some $c_{ij}\in \kp$. For a fixed $k$, there exists $i,j$ such that $mi+qj=k$ if and only $i\equiv k m^{-1}[q]$. The proof follows straightforwardly. 
\end{proof}

\begin{cor}\label{c} $
\kp((x))[y]=\apl \oplus y\apl \oplus\cdots \oplus y^{q-1}\apl=\apl \oplus x\apl \oplus\cdots \oplus x^{q-1}\apl.$ $\hfill\square$ \end{cor}

Notice that if $q>1$, neither $x$ nor $y$ belongs to $\apl$. 
Let us consider the union of all translated of $\apl$ and $\bpl$ by a monomial.
$$
\tapl=\bigcup_{i\in \zp,j\in \np}x^i y^j \apl,\qquad \tbpl=\bigcup_{i\in \zp,j\in \np}x^i y^j \bpl
$$
These sets are not stable by addition, but they both form a multiplicative monoid.

\begin{cor}\label{cor:size_support}
If $F\in \tbpl$ has bidegree $(d,n)$, its support has size $\Oc(d n/q)$.  $\hfill\square$
\end{cor}

In what follows, we simply say precision for Gauss precision. 

\subsubsection{Fast multiplication in $\tapl$.}
A key point for our purpose is that we have access to a faster multiplication  in $\tapl$ than in $\kp((x))[y]$. Let us start with an easy lemma. 

\begin{lem}\label{lem:mult}
Let $G,H\in \kp((x))[y]$ and let $N\in \zp$. The product $GH\mod x^N$ only depends on $G\mod x^{N-v_0(H)}$ and $H \mod x^{N-v_0(G)}$. 
\end{lem}

\begin{proof}
Clear.
\end{proof}

\begin{prop}\label{prop:mult}
Let $G,H\in \tapl$ of degree at most $d$. Given $n> 0$, we can compute $F=GH$ with precision $n+ v_0(F)$ with $\Ot(d n/q)$ operations in $\kp$. 
\end{prop}

\begin{proof}
Thanks to the relation $v_0(F)=v_0(G)+v_0(H)$, Lemma \ref{lem:mult} shows that it's enough to compute $F_0=G_0H_0$ where 
$$
G_0=G \mod x^{n+v_0(G)},\quad H_0=H \mod x^{n+v_0(G)}.
$$
Since $G_0,H_0\in \tbpl$, the supports of $G_0,H_0$ have size $\Oc(d n/q)$. Since $\tbpl$ is a monoid, the support of $F_0=G_0 H_0 \in \tbpl$ has also size $\Oc(d n/q)$. It follows from \cite[Proposition 6]{HoLe13} or \cite[Theorem 12]{HoLeSch13} that $F_0$ can be computed in time $\Ot(d n/q)$.
\end{proof}

We thus gain a factor $q$ when compared to usual bivariate multiplication. Note that fast multiplication of polynomials with prescribed support is based on a sparse multivariate evaluation-interpolation strategy (see \cite{HoLeSch13,HoLe13} and references therein), the crucial point here being that $F=GH$ remains sparse thanks to the monoid structure of $\tapl$. 

\subsubsection{Fast division in $\apl$.} Since the map $\tau_\lambda$ preserves the degree in $y$, both rings $\apl,\apl^+$ are euclidean rings when considering division with respect to $y$. Namely, given $F,G\in \kp((x))[y]$, the euclidean division $F=QG+R$, $\deg(R)<\deg(G)$ forces the euclidean division of $\hat{F},\hat{G}$ defined by \eqref{eq:tau_lambda} to be
$$
\hat{F}=\hat{Q}\hat{G}+\hat{R},\quad \hat{Q},\hat{R}\in \apl,\quad \deg(\hat{R})<\deg(\hat{G}).
$$
Moreover, the next lemma ensures that if $\hat{F},\hat{G} \in \apl^+$ (resp. $\bpl$) with $\hat{G}$ monic, then $\hat{Q},\hat{R} \in \apl^+$ (resp. $\bpl$).  

\begin{lem}\label{lem:prec_division} Let $F,G\in \kp((x))[y]$ with euclidean division $F=QG+R$. Assume that the leading coefficient of $G$ has valuation $v_0(G)$. Then 
$$
v_0(Q)\ge v_0(F)-v_0(G) \quad {\rm and} \quad v_0(R)\ge v_0(F).
$$
\end{lem}

\begin{proof}
See e.g. \cite{PoWe22} (a similar result holds for an arbitrary  valuation).
\end{proof}

Given $F\in \kp((x))[y]$ of degree $d$, let us denote by $\tilde{F}=y^{d} F(x,y^{-1})$ its reciprocal polynomial. We will need the following lemma.

\begin{lem}\label{lem:reciprocal} Let $F\in \apl$ of degree $d$. Then $\tilde{F}\in y^r\mathbb{A}_{-\lambda}$ where $r=d\mod q$. 
\end{lem}

\begin{proof}
Let $F\in \apl$ with expression as in Lemma \ref{lem:the_ring_apl}. Then 
$$
\tilde{F}=\sum_{k} \tilde{f_k}(y^q) y^{d-q\deg(f_k)-\alpha_\lambda(k)}  x^k,\quad f_k\in \kp[y]
$$
We have $d-q\deg(f_k)-\alpha_\lambda(k)\equiv  r+\alpha_{-\lambda}(k)\mod q$ and the claim follows from Lemma \ref{lem:the_ring_apl} applied in the ring $\mathbb{A}_{-\lambda}$.
\end{proof}

\begin{prop}\label{prop:div}
Let $F,G\in \apl$ of degree at most $d$, and suppose that the leading coefficient of $G$ has valuation $v_0(G)$. Given $n\ge 0$, we can compute $Q,R\in \apl$ with $\deg(Q)<\deg(G)$ such that 
$$F=QG+R \mod x^{v_0(F)+n}$$
with $\Ot(d n/q)$ operations in $\kp$.  
\end{prop}

\begin{proof}
Let $e=\deg(G)$ and $d=\deg(F)$. Assume $d>e$. Let us first reduce to the case $G$ monic.
We need to take care that multiplication by an arbitrary power of $x$ is not allowed in $\apl$. We proceed as follows. Let $k=-v_0(G)$ and let $\alpha=\alpha_{\lambda}(k)$. By Lemma \ref{lem:the_ring_apl}, we have $x^{k} y^\alpha\in \apl$ so the polynomials
$$
G_0=x^k y^\alpha G\quad {\rm and}\quad F_0=x^k y^\alpha F
$$
belong to $\apl$, with now $v_0(G_0)=0$. We are reduced to solve 
$$
F_0=Q G_0 +R_0 \mod x^{v_0(F_0)+n}
$$
in $\apl$, recovering $R$ for free from the relation $R_0=x^k y^\alpha R$. 
By assumption the leading coefficient $u(x)$ of $G_0$ is invertible in $\kp[[x]]$. Moreover, $\deg(G_0)=e+\alpha$ is divisible by $q$ and it follows  from Lemma \ref{lem:the_ring_apl} that $u\in \kp[[x^q]]\subset \apl$. Hence so does $u^{-1}$. Thus $u$ can be invert in $\apl^+$ with precision $n$ in time $\Ot(n/q)$, and we may suppose safely that $G_0$ is monic. Note that $$\deg(F_0)-\deg(G_0)=\deg(F)-\deg(G)=d-e.$$
The classical fast euclidean division $F_0=Q G_0 +R_0$ runs as follows:
\begin{enumerate}
\item Truncate $F_0$ at precision $n+v_0(F_0)$ and $G_0$ at precision $n$. 
\item Compute $\tilde{H}_0:=\tilde{G}_0^{-1}\mod y^{d-e+1}$ with precision 
$n$. 
\item Compute $\tilde{Q}=\tilde{F}_0 \tilde{H}_0\mod y^{d-e+1}$  with precision $n+v_0(F)$.
\item Compute $Q=y^{d-e}\tilde{Q}(x,y^{-1})$.
\item Compute $R_0=F_0-Q G_0$ with precision $n+v_0(F)$.
\end{enumerate}
Note that step 2 makes sense: since $G_0$ is monic, we have  $\tilde{G}_0(0)=1$ so $\tilde{G}_0$ can be invert in $\kp[[x]][[y]]$. This algorithm returns the correct output $F=QG+R$ if we do not truncate, see e.g. \cite[Thm 9.6]{GaGe13} and Lemma \ref{lem:prec_division} and Lemma \ref{lem:mult} ensure that truncations are correct to get $F_0=QG_0+R_0\mod x^{n+v_0(F)}$. Using quadratic Newton iteration, the inversion of $\tilde{G}_0$ at step 2 requires $\Oc(\log(d))$ multiplications and additions in $\kp[[x]][y]$ of degrees at most $d-e$ with  precision $n$ (see e.g. \cite[Thm 9.4]{GaGe13}). Since $q$ divides $\deg(G_0)$, Lemma \ref{lem:reciprocal} gives $\tilde{G}_0\in \mathbb{A}_{-\lambda}$, which is a ring. Hence all additions and multiplications required by \cite[Algorithm 9.3]{GaGe13} take place in $\ap_{-\lambda}$  and the cost of step 2 fits in the aimed bound thanks to Proposition \ref{prop:mult}. Since $\tilde{H}_0,\tilde{F}_0\in \tapl$ by Lemma \ref{lem:reciprocal}, we compute $\tilde{Q}$ at step 3 in time $\Ot((d-e)n/q)$ by Proposition \ref{prop:mult}. Step 4 is free. At step 5, the equation has degree $d+\alpha$ and vanish $\mod y^\alpha$, so its sparse size is $\Ot(dn/q)$ and step 5 fits too in the aimed bound since $F_0,Q,G_0\in \apl$. 
\end{proof}

%

\subsubsection{Fast Hensel lifting in $\apl$.} 

\begin{prop}\label{prop:hensel_apl}
Let $F\in \apl^+$ of degree $d$ and consider a coprime factorization $F(0,y)=f_0\cdots f_r f_\infty$ in $\apl\cap\kp[y]=\kp[y^q]$ with $f_i$ monic and $f_\infty=c\in \kp^\times$. Then there exists uniquely determined polynomials $F_0,\ldots,F_r,F_\infty \in \apl^+$ such that
$$
F=F_0\cdots F_r F_\infty,\quad  F_i(0,y)=f_i(0,y)\quad i=0,\ldots,k,\infty
$$
with $F_i$ monic of degree $\deg(f_i)$. We can compute the $F_i$'s with precision $n$ within $\Ot(dn/q)$ operations in $\kp$. Moreover, the truncated polynomials $F_i\mod x^n$ are uniquely determined by the equality $F\equiv F_0\cdots F_r F_\infty \mod x^n$. 
\end{prop}

\begin{proof}
This is the classical fast multi-factor hensel lifting, see e.g. \cite[Algorithm 15.17]{GaGe13}. The algorithm is based on multiplications and divisions of polynomials at precision $n$. The initial Bezout relations holds here in $\kp[y^q]\subset \apl$, and it follows that at each Hensel step, the input polynomials belong to the ring $\apl$. Moreover, all euclidean divisions satisfy the hypothesis of Proposition \ref{prop:div}. The claim thus follows from Proposition \ref{prop:mult} and Proposition \ref{prop:div} together with \cite[Theorem 15.18]{GaGe13}. Unicity of the lifting mod $x^n$ follows from \cite[Theorem 15.14]{GaGe13}.
\end{proof}

\begin{rem}
It is crucial to consider the factorization of $F(0,y)$ in the ring $\apl$. Typically, a polynomial of shape $y^q-1$ should be considered irreducible. Otherwise, the complexity will be $\Ot(dn)$ due to the loss of sparse arithmetic.
\end{rem}

\begin{rem}
Propositions \ref{prop:div} and \ref{prop:hensel_apl} appear also in \cite[Propositions 11 and 12]{PoRy15} under the assumption that $F\in bpl$ is monic. However, the proofs have not been published up to our knowledge. 
\end{rem}

\subsection{Fast $\vl$-adic Hensel lemma}

By the isomorphism $\tau_\lambda:\kp((x))[y]\to \apl$, the previous results translate in an obvious way in quasi-linear complexity estimates for $\vl$-adic truncated multiplication and division in $\kp((x))[y]$.

\begin{cor}
Let $\lambda\in \qp$ and let $G,H\in \kp((x))[y]$ of degree at most $d$. We can compute $F=GH$ at $\lambda$-precision $\vl(F)+\sigma$  with $\Ot(d\sigma)$ operations in $\kp$. 
\end{cor}

\begin{proof}
Follows from \eqref{eq:v0vl} together with Proposition \ref{prop:mult}.
\end{proof}

\begin{cor}\label{cor:lambda_mult}
We can multiply arbitrary polynomials $G,H\in \kp[x,y]$ in quasi-linear time with respect to the $\lambda$-size of the output. $\hfill\square$
\end{cor}

\begin{rem}
We could have used directly a sparse multivariate evaluation-interpolation strategy on the input polynomials $G,H$. However, we believe that using fast arithmetic in the ring $\apl$ is more convenient and offers more applications. 
\end{rem}

\begin{defi}\label{def:lambda_monic}
We say that $G\in \kp((x))[y]$ is $\lambda$-monic if its leading monomial $u y^e$ satisfies $\vl(uy^e)=\vl(G)$. 
\end{defi}

\begin{prop}\label{prop:lambda_div}
Let $F,G\in \kp((x))[y]$ of degrees at most $d$ with $G$ $\lambda$-monic. We can compute $Q,R\in \kp((x))[y]$ such that
$$
\vl(F-(QG+R))\ge \vl(F)+\sigma
$$
within $\Ot(d\sigma)$ operations. 
\end{prop}

\begin{proof}
We apply Proposition \ref{prop:div} to the polynomials $\hat{F}=\tl(F)$ and $\hat{G}=\tl(G)$. We are reduced to compute $\hat{Q},\hat{R}$ such that
\begin{equation}\label{eq:div_euclidienne_lambda}
v_0(\hat{F}-(\hat{Q}\hat{G}+\hat{R}))\ge v_0(\hat{F})+q\sigma.
\end{equation}
Since $G$ is assumed to be $\lambda$-monic, the leading coefficient $u(x)$ of $\hat{G}$ has valuation 
$v_0(\hat{G})$ and we conclude thanks to Proposition \ref{prop:div}.
\end{proof}

We get finally a fast Hensel lifting with respect to the valuation $\vl$.

\begin{defi}
Let $F=\sum c_{ij} x^i y^j\in \kp((x))[y]$ and $\sigma\in \frac{1}q \zp$. The $\lambda$-\emph{homogeneous component} of $F$ of degree $\sigma$ is
$$
F_\sigma=\sum_{i+j\lambda = \sigma} c_{ij} x^i y^j \in \kp[x^{\pm 1}][y].
$$
The $\lambda$-\emph{initial part} of $F$ is the $\lambda$-homogeneous component of $F$ of lowest degree $\vl(F)$, denoted by $\In_\lambda(F)$.  
\end{defi}

Let $F\in \kp((x))[y]$. The irreducible factorization of the $\lambda$-initial part of $F$ can be written in  a unique way (up to permutation) as
\begin{equation}\label{eq:init_facto}
\inl(F)= p_0 p_1\cdots p_k p_\infty \in \kp[x^{\pm 1}][y]
\end{equation}
where  $p_0=y^n$ with $n\in \np$, $p_\infty=u x^a$ with $a\in \zp$, $u\in \kp^\times$ and where $p_1,\ldots,p_k$ are coprime powers of irreducible $\lambda$-homogeneous monic polynomials, not divisible by $y$. 
The following result is well known (see e.g. \cite[Chapter VI]{Ca80}).

\begin{prop}\label{prop:lambda-hensel}
There exists unique polynomials $P_i^*\in \kp((x))[y]$ such that 
$$
F=P_0^*\cdots P_k^* P_\infty^* \in \kp((x))[y],\qquad \inl(P_i^*)=p_i,
$$
with $P_i^*$ monic of $\deg(P_i^*)=\deg(p_i)$ for $i=0,\ldots,k$. Moreover, the polynomial $P_i^*$ is $\lambda$-monic, and irreducible if $p_i$ is irreducible for all $i=0,\ldots,k$.
\end{prop}

\medskip
We get the following complexity result:

\begin{prop}\label{prop:approx_Gi}
Given $\sigma\in \qp^+$, and given the irreducible factorization \eqref{eq:init_facto} we can compute $P_0,\ldots,P_\infty\in \kp[x^{\pm 1}][y]$ such that 
$$
\vl(P_i^*-P_i)> \vl(P_i^*)+\sigma\quad \forall\,i=0,\ldots,k,\infty.
$$
in time $\Ot(d\sigma)$. We have then
$$
\vl(F-P_0\cdots P_\infty)> \vl(F)+\sigma.
$$
\end{prop}

\begin{proof}
Up to multiply $F$ by a suitable monomial $x^i y^\alpha$ with $0\le \alpha<q$, we may assume that $\vl(F)=0$. Then we apply Proposition \ref{prop:hensel_apl} to the polynomial $\hat{F}=\tau_\lambda(F)$, starting from the factorization of $\hat{F}(0,y)$ induced by the factorization of $\inl(F)$, and with a suitable Gauss precision in order to recover the desired $\lambda$-precision. The cost and the unicity of the truncated polynomial $P_i$ follow from Proposition \ref{prop:hensel_apl}. 
\end{proof}

\begin{rem}
This result improves \cite[Corollary 2]{PoWe22} which gives the complexity estimate $\Ot(d(\sigma+\vl(F))$. Proposition \ref{prop:hensel_apl} has a significant impact when needed a lifting precision closed to $\vl(F)$. This is precisely the case for our application to bivariate factorization. 
\end{rem}

\begin{defi}\label{def:partial_facto}
We denote \emph{\texttt{PartialFacto($F,\lambda,\sigma$)}} the algorithm which computes the factorization \eqref{eq:init_facto} of $\inl(F)$ and returns the truncated factors $P_0,\ldots,P_\infty$ following Proposition \ref{prop:approx_Gi}.
\end{defi}

\subsection{Fast $\vl$-adic factorization}\label{ssec:fast_lambda_facto}

We want now to compute the complete irreducible factorization of $F$ in $\kp((x))[y]$. Although our target precision is measured in terms of the valuation $\vl$, we will perform recursive calls of \texttt{PartialFacto} with various valuations $\vlp$. The integer $\ml(F)$ introduced in \eqref{eq:mlambda_intro} will play a key role. 

\subsubsection{The $\lambda$-defect of straightness}\label{ssec:defect_of_straightness}

%
%
%
%

\begin{defi}\label{def:mlambda}
Given $P=\sum_{i=s}^n p_i y^i\in \kp((x))[y]$ with $p_s,p_n\ne 0$, we denote $\iny(P)=p_s y^s$ the initial term of $P$ and $\lty(P)=p_n y^n$  the leading term of $P$.  We define 
$$\al(P)=\vl(\iny(P))-\vl(P)\quad {\rm and}\quad \bl(P)=\vl(\lty(P))-\vl(P).$$
The \emph{$\lambda$-defect of straightness} of $P$ is $\ml(P)=\max(\al(P),\bl(P)).$
\end{defi}

Recall from the introduction that $\Lambda(P)$ is the lower convex hull of the set of points $(i,v_0(p_i))$, $i=s,\ldots,n$, where $v_0(p_i)$ is the $x$-adic valuation. By convexity, $v_0(p_i)+i\lambda$ takes its maximal value at $i=s$ or $i=n$, hence the definition of $\ml(P)$ coincides with \eqref{eq:mlambda_intro}. The terminology for $\ml$ is justified by the following fact:

\begin{lem}\label{lem:propriety_mlambda}
The following properties hold:
\begin{enumerate}
\item  $\al(P)\ge 0$.
\item $\bl(P)\ge 0$ with equality if and only if $P$ is $\lambda$-monic.
\item $\ml(P)\ge 0$  with equality if and only if $\Lambda(P)$ is one-sided of slope $-\lambda$. 
\end{enumerate}
\end{lem}

\begin{proof}
This follows from the equality
$\vl(P)=\min\{\vl(p_i y^i),\, i=s,\ldots,n\}$.
\end{proof}

\begin{cor}\label{cor:mlPQ}
Let $P,Q\in \kp((x))[y]$. We have $\Lambda(PQ)=\Lambda(P)+\Lambda(Q)$ and
$$\ml(PQ)\ge \max(\ml(P),\ml(Q))$$ with equality if $\Lambda(Q)$ or $\Lambda(P)$ is one-sided of slope $-\lambda$.
\end{cor}

\begin{proof}
First equality is a well known variant of Ostrowski's theorem. Since $\iny$ and $\lty$ are multiplicative operators and $\vl$ is a valuation, we get $$\al(PQ)=\al(P)+\al(Q)\quad{\rm and}\quad \bl(PQ)=\bl(P)+\bl(Q).
$$ The inequality for $\ml$ follows straightforwardly. If $\Lambda(Q)$ or $\Lambda(P)$ is one-sided of slope $-\lambda$, the equality follows from point (3) of Lemma \ref{lem:propriety_mlambda}.
\end{proof}

\subsubsection{Comparisons between various valuations}

\begin{lem}\label{lem:vl_vs_vlp}
Let $\lambda'\ge \lambda$ and let $P\in\kp((x))[y]$ of degree $n$. Then 
$$
\vl(P)\le \vlp(P) \le \vl(P)+n(\lambda'-\lambda)
$$
\end{lem}

\begin{proof}
Since $i+j\lambda\le i+j\lambda'$ we get immediately $\vl\le \vlp$. 
Let $(i_0,j_0)$ in the support of $P$ such that $\vl(P)=i_0+j_0\lambda$. We get
$$
\vlp(P)\le i_0+j_0\lambda' = i_0+j_0\lambda + j_0(\lambda'-\lambda) = \vl(P)+j_0(\lambda'-\lambda)
$$
and we conclude thanks to $j_0\le n$.
\end{proof}

\begin{defi}
Let $P_0,P\in\kp((x))[y]$. We say that $P_0$ approximate $P$ with relative $\lambda$-precision $\sigma$ if $\vl(P-P_0)>\vl(P)+\sigma$. We say that $P$ is known with relative $\lambda$-precision $\sigma$ if we know such an approximant $P_0$.
\end{defi}

\begin{cor}\label{cor:precprime}
Let $P\in\kp((x))[y]$ of degree $n$, let $\lambda,\lambda'\in\qp$ and let $\sigma\ge 0$.
If $P$ is known with relative $\lambda'$-precision 
\begin{equation}\label{eq:sigmaprime}
\sigma'=\sigma'(\lambda,\lambda',\sigma,P):=\begin{cases} \sigma+\vl(P)-\vlp(P)+ n(\lambda'-\lambda) \quad{\rm if}\,\, \lambda'\ge \lambda \\
\sigma+\vl(P)-\vlp(P)  \qquad\,\qquad\qquad{\rm if}\,\, \lambda'\le \lambda 
\end{cases}
\end{equation}
then $P$ is known with relative $\lambda$-precision $\sigma$.
\end{cor}

\begin{proof}
The first claim follows from the second inequality in Lemma \ref{lem:vl_vs_vlp} for the case $\lambda'\ge \lambda$ and from the first inequality in Lemma \ref{lem:vl_vs_vlp} for the case $\lambda'\le \lambda$. 
\end{proof}

\begin{lem}\label{lem:lm_vs_lmprime}
We keep notations of Corollary \ref{cor:precprime}. 

$\bullet$ If $\lambda'\ge \lambda$, then 
$
\ml(P)+\vl(P)-n\lambda \ge \mlp(P)+\vlp(P)-n\lambda'.
$

$\bullet$ If $\lambda'\le \lambda$, then 
$
\ml(P)+\vl(P) \ge \mlp(P)+\vlp(P).
$
\end{lem}

\begin{proof}
Denoting $P=\sum_{i=s}^n p_i y^i$, the first inequality is equivalent to that 
$$
\max(v_0(p_s)-(n-s)\lambda,v_0(p_n))\ge \max(v_0(p_s)-(n-s)\lambda',v_0(p_n)),
$$
which follows from the assumption $\lambda\le \lambda'$. The second inequality is equivalent to 
$$
\max(v_0(p_s)+s\lambda,v_0(p_n)+n\lambda)\ge \max(v_0(p_s)+s\lambda',v_0(p_n)+n\lambda'),
$$
which follows from the assumption $\lambda\ge\lambda'$.
\end{proof}

\begin{cor}\label{cor:sigmaprime_vs_mprime}
We have $\sigma'-\sigma\ge \mlp(P)-\ml(P).$
\end{cor}

\begin{proof}
Combining \eqref{eq:sigmaprime} and Lemma \ref{lem:lm_vs_lmprime} leads to the desired inequality.
\end{proof}

We will need also an upper bound for $\sigma'$ in terms of $\sigma$. 

\begin{lem}\label{lem:bound_for_sigmaprime}
We keep notations of Corollary \ref{cor:precprime}. 
Suppose that $\bl(P)=0$ if $\lambda'\ge \lambda$ and that $\al(P)=0$ and $P$ not divisible by $y$ if $\lambda'\le\lambda$. Then
$
\sigma'\le \sigma + \mlp(P).
$
\end{lem}

\begin{proof}
Suppose $\lambda'\ge \lambda$. By \eqref{eq:sigmaprime}, the inequality is equivalent to 
$$
\mlp(P)\ge \vl(P)-\vlp(P)+n(\lambda'-\lambda)
$$
Both sides are invariant when multiplying $P$ by a power of $x$, hence we can safely suppose $P$ monic in $y$. The hypothesis $\bl(P)=0$ is still true and we get $\vl(P)=\vl(y^n)=n\lambda$. We are reduced to show 
that
$
\mlp(P)\ge n\lambda'-\vlp(P)=\blp(P)$, which 
follows from Definition \ref{def:mlambda}.  Suppose now $\lambda'\le \lambda$. By \eqref{eq:sigmaprime}, we need to show that
$
\mlp(P)\ge \vl(P)-\vlp(P) .
 $ 
By hypothesis, we have $\vl(P)=\vl(p_0)=\vlp(p_0)$. We are reduced to show that $\mlp(P)\ge \vlp(p_0)-\vlp(P)=\alp(P)$, which follows from Definition \ref{def:mlambda}.
\end{proof}

\subsubsection{Recursive calls}

Let $F\in \kp((x))[y]$. We fix $\lambda\in \qp$ and a relative $\lambda$-precision $\sigma\ge 0$. Following Definition \ref{def:partial_facto}, let us consider 
$$L=[P_0,P_1,\ldots,P_k,P_\infty]=\,\textrm{\texttt{PartialFacto}}(F,\lambda,\sigma). 
$$
Assuming $F$ non degenerated, we know thanks to Proposition \ref{prop:lambda-hensel} that $P_0,\ldots, P_\infty$ approximate some coprime factors $P_0^*,P_1^*,\ldots,P_k^*,P_\infty^*$ of $F$ with relative $\lambda$-precision $\sigma$. Moreover, the polynomials $P_1^*,\ldots,P_k^*$ and their approximant are 
 irreducible. There remains to factorize (if required) the polynomials $P_0^*$ and $P_\infty^*$. We denote for short
$$
(G^*,H^*)=(P_0^*,P_\infty^*)\quad {\rm and}\quad (G,H)=(P_0,P_\infty).
$$

\begin{lem}\label{lem:mlP}
The polynomials $G$ and $G^*$ are monic of same degree and $H$ and $H^*$ are not divisible by $y$. 
Moreover:

$\bullet$ If $P$ divides $G$, then $\ml(P)=\al(P)$ and $\bl(P)=0$.

$\bullet$ If $P$ divides $H$, then $ \ml(P)=\bl(P)$ and $\al(P)=0$.
\end{lem}

\begin{proof} The first claim follows from Proposition \ref{prop:approx_Gi}, Proposition \ref{prop:lambda-hensel} and \eqref{eq:init_facto}. More precisely, denoting  $G=c_0+\cdots + c_n y^n$ and $H=h_0+\cdots + h_m y^m$ with $c_n,h_m\ne 0$, we have
$$
\inl(G)=\inl(G^*)=\inl(y^n)\quad {\rm and}\quad \inl(H)= \inl(H^*)=\inl(h_0).
$$
As $\vl(G-\inl(G))>\vl(G)$ we deduce $\bl(G)=0$ and $\ml(G)=\al(G)$. In the same way, we get $\al(H)=0$ and $\ml(H)=\bl(H)$. 
If $P$ divides $G$, we have $\bl(P)\le \bl(G)$ by multiplicativity of $\bl$. As $\bl(P)\ge 0$, this forces $\bl(P)=0$, and thus $\ml(P)=\al(P)$. 
If $P$ divides $H$, then $0\le \al(P)\le \al(H)$ forces now $\al(P)=0$ and  $\ml(P)=\bl(P)$. 
\end{proof}

We need a lower bound on $\sigma$ which ensures that we can detect the irreducible factors of $G^*$ and $H^*$ on their approximants $G$ and $H$. 


\begin{lem}\label{lem:GvsGstar}
Suppose that $\sigma\ge \ml(F)$. Then $\Lambda(G)=\Lambda(G^*)$ and the restriction of $G$ and $G^*$ to their lower convex hull coincide. The same assertion is true for $H$ and $H^*$. In particular, $\ml(G)=\ml(G^*)$ and $\ml(H)=\ml(H^*)$.
\end{lem}

\begin{proof}
By Lemma \ref{lem:mlP}, we have  $G^*=c_s^* y^s+\cdots+ y^n$ with $c_s^*\ne 0$ and $G=c_s y^s+\cdots +  y^n$ (we might have \textit{a priori} $c_s=0$). By a convexity argument, we are reduced to show that  $c_s,c_s^*\in \kp((x))$ have same $x$-adic initial term. We have
$$\vl(c_s^*y^s -c_s y^s)\ge \vl(G-G^*)> \vl(G^*)+\sigma\ge \vl(G^*)+\ml(G^*)=\vl(c_s^*y^s),
$$
the first inequality by definition of $\vl$, the second inequality by Proposition \ref{prop:lambda-hensel}, the last inequality by hypothesis $\sigma\ge \ml(F)$ combined with Corollary \ref{cor:mlPQ}, and the last equality by 
Lemma \ref{lem:propriety_mlambda} since $G^*$ is $\lambda$-monic (Proposition \ref{prop:lambda-hensel}). We deduce that $\inl(c_s^*y^s)=\inl(c_s y^s)$, from which it follows that $c_s,c_s^*\in \kp((x))$ have same initial term as required. 
The assertion for $H$ is proved in the same way, focusing now on the leading term of $H$. 
\end{proof}

Assuming $F$ non degenerated, its irreducible factorization in $\kp((x))[y]$ is deduced from the irreducible factorization of its lower edges polynomials. Hence, Lemma \ref{lem:GvsGstar} ensures that knowing $G$ and $H$ at precision $\sigma\ge \ml(F)$ is sufficient to detect all remaining \emph{irreducible} factors of $F$.

\subsubsection{Divide and conquer}\label{sssec:divandconquer}
We apply now recursively  \texttt{PartialFacto} to $G$ and $H$  with respect to some well chosen slopes $\lambda_G$ and $\lambda_H$.

\begin{defi}\label{def:average_slope}
Let $n>s$. The \emph{average slope} of $P=\sum_{i=s}^{n} p_i y^i\in \kp((x))[y]$ is
$$
\lambda_P:=-\frac{v_0(p_n)-v_0(p_s)}{n-s} \in \qp.
$$
In other words, $-\lambda_P$ is the slope of the segment joining the two extremities of the lower boundary $\Lambda(P)$.  
\end{defi}

This slope is chosen so that the $\lambda_P$-valuation of the leading term and the initial term of $P$ coincide. Equivalently, it satisfies
\begin{equation}
\label{eq:mlambdaP=alambdaP}
m_{\lambda_P}(P)=a_{\lambda_P}(P)=b_{\lambda_P}(P).
\end{equation}
We deduce:

\begin{prop}\label{prop:lambdaG<lambda}
Let $\lambda$ and $G,H$ as above, and suppose $G,H$ of positive $y$-degree. 

$\bullet$ If $P$ divides $G$ then $\lambda_P\ge \lambda$.

$\bullet$ If $P$ divides $H$ then $\lambda_P\le \lambda$.

\noindent 
In both cases, we have $\mlP(P)\le \ml(P)$.
\end{prop}

\begin{proof}
We can write $P=a_0+\cdots +  a_n y^n$ with $a_0, a_n\ne 0$ and $n\ge 1$. If $P$ divides $G$, Lemma \ref{lem:mlP} implies $\vl(a_0)\ge \vl(a_n y^n)$. Thus $v_0(a_0)\ge v_0(a_n)+n\lambda$, which implies $\lambda_P\ge \lambda$. We get
$$
\mlP(P) = a_{\lambda_P}(P)=v_0(a_0)-v_{\lambda_P}(P) \le  v_0(a_0)-\vl(P) =\al(P) = \ml(P),
$$
the first equality by \eqref{eq:mlambdaP=alambdaP}, the inequality thanks to $ \lambda_P \ge \lambda$ and the last two equalities by Lemma \ref{lem:mlP}. If $P$ divides $H$,
 Lemma \ref{lem:mlP} forces now $\vl(P)=\vl(a_0)\le \vl(a_n y^n)$ and thus $\lambda_P\le \lambda$.  By \eqref{eq:mlambdaP=alambdaP}, we get
$$
\mlP(P) = b_{\lambda_P}(P) = v_{\lambda_P}(a_n y^n)-v_{\lambda_P}(P).
$$
On one hand, we have $v_{\lambda_P}(P)\ge v_0(a_0) = \vl(P)$. On the other hand $\lambda_P\le \lambda$ implies $v_{\lambda_P}(a_n y^n)\le v_{\lambda}(a_n y^n)$. We get
$$
\mlP(P) \le v_{\lambda}(a_n y^n)- \vl(P)=\bl(P) = \ml(P),$$
the two equalities by Lemma \ref{lem:mlP}.
\end{proof}

\begin{defi}\label{def:sigmaG}
Let $\lambda$ be fixed and let $P\in \kp((x))[y]$.  Given a $\lambda$-precision $\sigma$, we denote $\sigma_P=\sigma'(\lambda,\lambda_P,\sigma,P)$ the  precision induced by \eqref{eq:sigmaprime} with $\lambda'=\lambda_P$. 
\end{defi}

We deduce the following key uniform upper bound for $\sigma_P$. 

\begin{prop}\label{prop:key_mlG}
Suppose that $\sigma\ge \ml(F)$. If $P$ divides $G$ or $H$, then 
$$\mlP(P)\le \sigma_P \le\sigma+\ml(F).$$
\end{prop}

\begin{proof}
The inequality $\mlP(P)\le \sigma_P$ follows from Corollary \ref{cor:sigmaprime_vs_mprime}. By Lemma  \ref{lem:bound_for_sigmaprime}, we get  $\sigma_P \le\sigma+\mlP(P)$. By Proposition  \ref{prop:lambdaG<lambda}, we get $\mlP(P)\le \ml(P)$. Since $P$ divides $G$, we have $\ml(P)\le \ml(G)$ by Corollary \ref{cor:mlPQ}. By Lemma \ref{lem:GvsGstar}, we have $\ml(G)=\ml(G^*)$, and Corollary \ref{cor:mlPQ} again gives $\ml(G^*)\le \ml(F)$. The claim follows.
\end{proof}

The last key result ensures that using the slopes $\lambda_G$ and $\lambda_H$ lead to a divide and conquer strategy. 
Given $P\in \kp((x))[y]$, we denote in what follows by $V_P$ the euclidean volume of the convex hull of $\Lambda(P)$.

\begin{prop}\label{prop:DGDH}
Let $F \in \kp((x))[y]$ and suppose that $\lambda=\lambda_F$ (as it will be the case at the recursive calls). Let $G,H\in \kp((x))[y]$ as defined above.  
\begin{enumerate}
\item $d \,m_{\lambda_F}(F)/2 \le V_F\le d \,m_{\lambda_F}(F)$. 
\item $V_F=0$ if and only if $\Lambda(F)$ is one-sided, in which case its slope is $\lambda_F$. 
\item We have $(V_G+V_H)\le V_F/2$.
\end{enumerate}
\end{prop}

\begin{proof}
We still denote $\Lambda(F)$ the convex hull of the lower boundary $\Lambda(F)$. Let $ABCD$ be the smallest parallelogram with two vertical sides containing $\Lambda(F)$ such that $C$ and $D$ are respectively the left end point and the right end point of $\Lambda(F)$ and $(AD)$ and $(BC)$ are vertical (figure \ref{volume} below).
\begin{figure}[h]
	\caption{Illustrated proof of Proposition \ref{prop:DGDH}.}
	\label{volume}
	\bigskip
	\scalebox{0.20}{\begin{picture}(0,0)%
\includegraphics{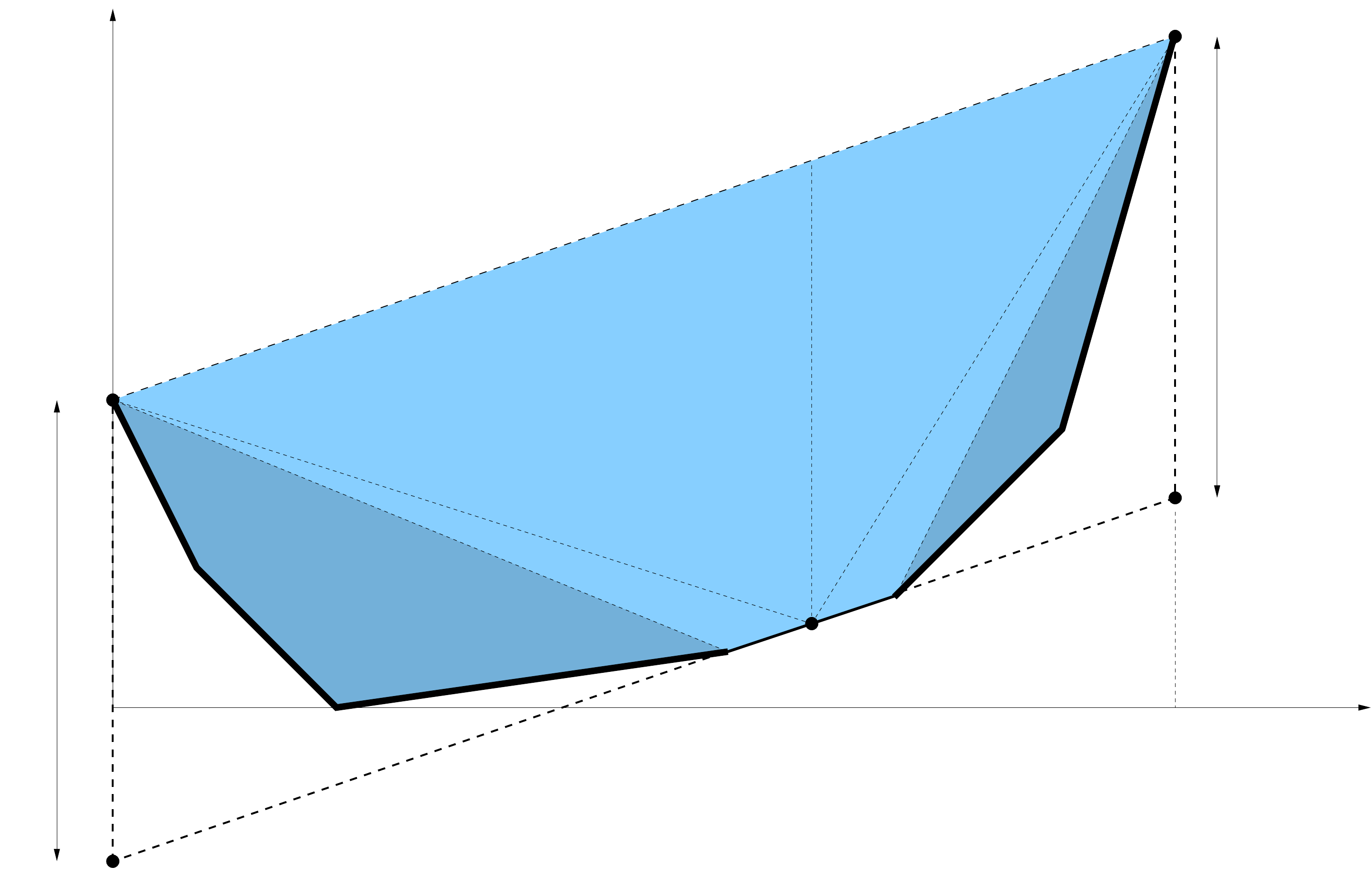}%
\end{picture}%
\setlength{\unitlength}{4144sp}%
\begingroup\makeatletter\ifx\SetFigFont\undefined%
\gdef\SetFigFont#1#2#3#4#5{%
  \reset@font\fontsize{#1}{#2pt}%
  \fontfamily{#3}\fontseries{#4}\fontshape{#5}%
  \selectfont}%
\fi\endgroup%
\begin{picture}(22077,14353)(886,-13841)
\put(5626,-9286){\scalebox{3.5}{$\Lambda(G)$}}
\put(16750,-6361){\scalebox{3.5}{$\Lambda(H)$}}
\put(11026,-6136){\scalebox{3.5}{$\Lambda(F)$}}
\put(19576,-11536){\scalebox{3.5}{$d=\deg(F)$}}
\put(2151,-13786){\scalebox{3.5}{$A$}}
\put(20026,-7936){\scalebox{3.5}{$B$}}
\put(20026,389){\scalebox{3.5}{$C$}}
\put(2026,-5911){\scalebox{3.5}{$D$}}
\put(13951,-9961){\scalebox{3.5}{$E$}}
\put(20926,-3886){\scalebox{3.5}{$b_{\lambda}(F)=m_{\lambda}(F)$}}
\put(201,-9961){\scalebox{3.5}{$a_{\lambda}(F)$}}
\end{picture}%
}
\end{figure}
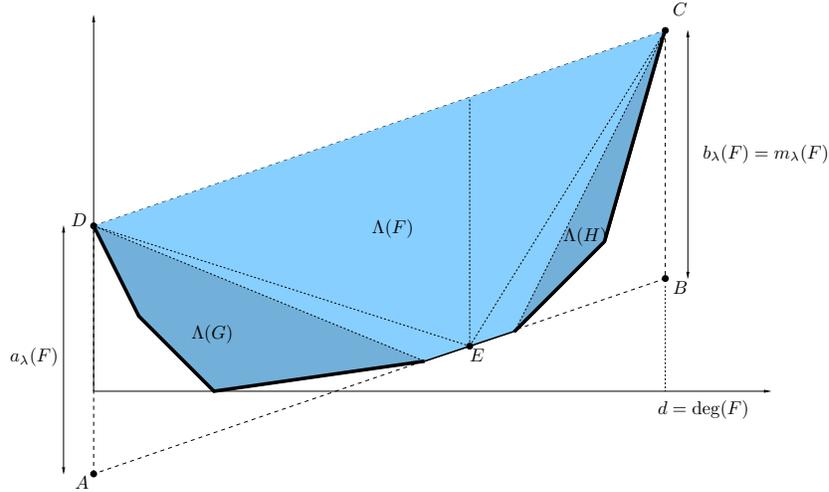

\noindent
Denote $\lambda=\lambda_F$ for short. The line $(AB)$ has equation $i+j\lambda=\vl(F)$, and the segments $[AD]$ and $[BC]$ have both length $\al(F)=\bl(F)=\ml(F)$ (Definition \ref{def:mlambda}) by choice of the average slope. Hence $ABCD$ has volume $d \ml(F)$, which gives $V_F\le d\,\ml(F)$. By construction, there exists $E\in [AB]\cap \Lambda(F)$ and by convexity, the triangle $CDE$ is contained in $\Lambda(F)$.  Since $\Vol(CDE)=\frac{1}2 \Vol(ABCD)$, the inequality $d \,m_{\lambda_F}(F)/2 \le V_F$ follows, proving first point. The second item is immediate. Since $\Lambda(G)$ is the Minkovski summand of $\Lambda(F)$ whose all minus slopes are strictly greater than $\lambda$ (Proposition \ref{prop:lambdaG<lambda}),  we may suppose that (up to translation) $\Lambda(G)\subset \Lambda(F)$ with left end point $C$ and right end point $I\in [AE]$. By convexity, $\Lambda(G)\subset AED$ and $V_G\le \Vol(AED)$. In the same way, we find $V_H\le \Vol(BCE)$. On the other hand, we have $V_F\ge V_G+V_H+\Vol(CDE)$. We conclude thanks to the relation $\Vol(CDE)=\Vol(AED)+\Vol(BCE)$.
\end{proof}

\begin{rem}
The partial factorization of $F$ with respect to $\lambda_F$ is $F=G^* Q^* H^*$ where $Q^*=P_1^*\cdots P_k^*$ has a one-sided lower boundary slope $\lambda_F$ (which is $[AB]\cap \Lambda(F)$ on figure \ref{volume}). However, although we use the terminology "slope", the rational $\lambda_F$ is generally not a slope of $\Lambda(F)$. In such a case, the intersection $[AB]\cap \Lambda(F)$ is reduced to a point and the partial $\lambda_F$-factorization of $F$ is simply $F=G^* H^*$. An important point is that $G^*$ and $H^*$ are not trivial factors as soon as  $\Lambda(F)$ has several slopes.
\end{rem}

\subsubsection{Proof of Theorem \ref{thm:lambda-facto}}

In the following algorithm, $\sigma_G,\sigma_H$ are defined by Definition \ref{def:sigmaG}, in terms of the input $(\lambda,\sigma)$ and the current slopes $\lambda_G$, $\lambda_H$. 

\begin{algorithm}[ht]
  \nonl\TitleOfAlgo{\texttt{Facto}($F, \lambda,\sigma$)\label{algo:divtype}}
  \KwIn{%
    $F\in \kp((x))[y]$ monic non degenerated, $\lambda\in \qp$ and $\sigma\ge \ml(F)$}%
  \KwOut{The irreducible factors $F$ with relative $\lambda$-precision $\sigma-\ml(F)$}%
 \lIf{$\deg(F)\le 1$}{\Return{$[F]$}}%
  $[P_0,P_1,\ldots,P_k,P_\infty]\gets$\texttt{PartialFacto}($F,\lambda,\sigma$)\label{step:partial_facto}\;%
  $G\gets P_0$, $H\gets P_\infty$\;%
    \leIf{$\deg(G)=0$}{$L_G\gets [\,\,]$}{$L_G\gets$\texttt{Facto}($G,\lambda_G,\sigma_G$)}\label{step:factoG}%
    \leIf{$\deg(H)=0$}{$L_H\gets [\,]$}{$L_H\gets$\texttt{Facto}($H,\lambda_H,\sigma_H$)}\label{step:factoH}%
    \Return{$[P_1,\ldots,P_k]\cup L_G\cup L_H$}%
\end{algorithm}

\begin{thm}\label{prop:fast_facto}
Given $F\in \kp((x))[y]$ non degenerate with irreducible factors $F_1^*,\ldots,F_s^*$ and given $\sigma\ge \ml(F)$, running \emph{\texttt{Facto}}($F,\lambda,\sigma$) returns a list of irreducible monic coprime polynomials $F_1,\ldots,F_s\in \kp[x,y]$ such that  
$$
\vl(F-F_1\cdots F_s)-\vl(F)> \sigma
$$
within $\Ot(d\sigma)$ operations in $\kp$.  Moreover,
$$
\vl(F_i-F_i^*)-\vl(F_i^*)> \sigma-\ml(F)
$$
for all $i=1,\ldots,s$.
\end{thm}

We will need the following lemma.

\begin{lem}\label{lem:approx}
If $\vl(A^*-A)>\vl(A)+\sigma$ and $\vl(B-B^*)>\vl(B)+\sigma$, then $\vl(A^* B^*-AB)>\vl(AB)+\sigma$.
\end{lem}

\begin{proof}
Follows from $A^* B^*-AB=A^*(B^*-B)+B(A^*-A)$ together with $\vl(A)=\vl(A^*)$ and $\vl(B)=\vl(B^*)$. 
\end{proof}

\medskip
\noindent
\textit{Proof of Theorem \ref{prop:fast_facto}.}

\smallskip
\noindent
$\bullet$ \textbf{Correctness.} By induction on the number of recursice calls. If the algorithm stops at step 2, then the result follows from Proposition \ref{prop:approx_Gi}. Else, we know that $G$ and $H$ are not degenerated (Lemma \ref{lem:GvsGstar}) and approximate $G^*$ and $H^*$ with relative $\lambda$-precision $\sigma$ (Proposition \ref{prop:approx_Gi}). As $\sigma_G\ge \mlG(G)$ (Proposition \ref{prop:key_mlG}), we deduce by induction that
\texttt{Facto}($G,\lambda_G,\sigma_G$) returns some approximants $G_1,\ldots,G_t$ of the irreducible factors $G_1^*,\ldots,G_t^*$ of $G$ such that
\begin{equation}\label{eq:eq1}
\vlG(G-G_1\cdots G_t)-\vlG(G)> \sigma_G,
\end{equation}
with moreover
\begin{equation}\label{eq:eq2}
\vlG(G_i-G_i^*)-\vlG(G_i^*)> \sigma_G-\mlG(G) \quad \forall \, i=1,\ldots,t.
\end{equation}
Corollary \ref{cor:precprime} and \eqref{eq:eq1} forces
$$
\vl(G-G_1\cdots G_t)-\vl(G)> \sigma.
$$
Since $\vl(G)=\vl(G^*)$ and $\vl(G-G^*)-\vl(G^*)> \sigma$, we deduce 
$$
\vl(G^*-G_1\cdots G_r)-\vl(G^*)> \sigma.
$$
In the same way, \texttt{Facto}($H,\lambda_H,\sigma_H$) computes an approximate irreducible factorization of $H^*$ such that
$$
\vl(H^*-H_1\cdots H_u)-\vl(H^*)> \sigma.
$$
We have $F=G^* P_1^*\cdots P_k^* H^*$ and we have too $\vl(P_i^*-P_i)-\vl(P_i^*)>\sigma$ (Proposition \ref{prop:approx_Gi}). The polynomials $(F_1,\ldots,F_s)=(P_1,\ldots,P_k,G_1,\ldots,G_t,H_1,\ldots,H_u)$ approximate the irreducible factors of $F$ and Lemma \ref{lem:approx} implies
$$
\vl(F-F_1\cdots F_r)-\vl(F)>\sigma
$$
as required. 
There remains to show that $\vl(F_i-F_i^*)-\vl(F_i^*)> \sigma-\ml(F)$ for all $i$. This is true for the factors $P_j$ by Proposition \ref{prop:lambda-hensel}. Let us consider a factor $A=G_i$. As $\lambda_G\ge \lambda$, \eqref{eq:eq2} combined with \eqref{eq:sigmaprime} gives
$$
\vlG(A-A^*)> \vlG(A)+\sigma+\vl(G)-\vlG(G)+d_G(\lambda_G-\lambda).
$$
Denote $B$ the (truncated) cofactor of $A$ in $G$. Using $\vl(A-A^*)+d_A(\lambda_G-\lambda)\ge \vlG(A-A^*)$ (Lemma \ref{lem:vl_vs_vlp}) together with $d_G=d_A+d_B$, $\vl(G)=\vl(A)+\vl(B)$ and $\vlG(G)=\vlG(A)+\vlG(B)$, the previous inequality implies that
$$
\vl(A-A^*)> \vl(A)+\sigma+\mlG(B)-\ml(B).
$$
As $\mlG(B)\ge 0$ and $\ml(B)\le \ml(F)$ (Corollary \ref{cor:mlPQ}), we get the desired inequality
$$
\vl(A-A^*)-\vl(A)>\sigma-\ml(F).
$$
 We prove in a similar way the analogous assertion if $A=H_i$ is a factor of $H$.

\smallskip
\noindent
$\bullet$ \textbf{Complexity.} There is at most $1+\lceil \log_2(V_F) \rceil =\Oc(\log_2(d\ml(F)))=\Oc(\log_2(d\sigma))$ recursive calls thanks to Proposition \ref{prop:DGDH} (the $+1$ due to the fact that the initial slope $\lambda$ is random).
At each level of the tree of recursive calls, the procedure \texttt{PartialFacto} is called on a set of polynomials $P$ dividing $G$ or $H$ and whose degree sum is at most $d_G+d_H\le d$, and with $\lambda_P$-precision $\sigma_P$ for each $P$. By  Proposition \ref{prop:key_mlG}, $\sigma_P\le \sigma+\ml(F)\le 2\sigma$ for all $P$, and we conclude thanks to Proposition \ref{prop:lambda-hensel}. $\hfill\square$
%
%
%

\bigskip
\noindent
\textit{Proof of Theorem  \ref{thm:lambda-facto}.}  Theorem \ref{thm:lambda-facto} follows straightforwardly from Theorem \ref{prop:fast_facto}, taking into account the cost of the factorizations \eqref{eq:init_facto} of the various quasi-homogeneous initial components. These factorizations are not trivial only when $\lambda$ is a slope of $\Lambda(F)$, in which case the degree of the underlying univariate factorization corresponds to the lattice length of the edge of slope $\lambda$. $\hfill\square$

\section{Application to convex-dense bivariate factorization}\label{sec:facto}


This section is dedicated to derive from Theorem \ref{thm:lambda-facto} a fast algorithm for factoring a bivariate polynomial $F\in \kp[x,y]$. 
 We follow closely \cite{We16}, which generalizes the usual factorization algorithm of  \cite{L2,L1} to the case $F(0,y)$ non separable.  To be consistent with \cite{L2,We16}, we denote from now on by $\F_i$ the factors of $F$ in $\kp((x))[y]$ and by $F_j$ the factors of $F$ in $\kp[x,y]$.


\subsection{The recombination problem}\label{ssec:recombination}

In all what follows, we assume that the input $F\in \kp[x,y]$ is primitive and separable of degree $d$ with respect to $y$ (see \cite{L3} for fast separable factorization). We normalize $F$ by requiring that its coefficient attached to the right end point of $\Lambda(F)$ equals $1$. Up to permutation, $F$ admits a unique factorization
\begin{equation}\label{factorat}
F=F_1\cdots F_\rho \in \kp[x,y],
\end{equation}
where each $F_j\in \kp[x,y]$ is irreducible and normalized. Also, $F$ admits a unique analytic factorization of shape
\begin{equation}\label{factoan}
F=u\F_1\cdots \F_s \in \kp[[x]][y],
\end{equation}
with $\F_i\in \kp[[x]][y]$ irreducible with leading coefficient $x^{n_i}$, $n_i\in \np$ and $u\in \kp[x]$, $u(0)\ne 0$. We thus have
\begin{equation}\label{prod}
F_j=c_j \F_1^{v_{j1}}\cdots \F_s^{v_{js}}, \,\quad j=1,\ldots,\rho,
\end{equation}
for some unique $v_{ji}\in\{0,1\}$, and with $c_j\in \kp[x]$, $c_j(0)=1$. 
The recombination problem consists to compute the exponent vectors 
$$
v_j=(v_{j1},\ldots,v_{js})\in \{0,1\}^{s}
$$
for all $j=1,\ldots,\rho$. Then, the computation of the $F_j$'s follows easily. Since $F$ is separable by hypothesis, the vectors $v_j$ form a partition of $(1,\ldots,1)$ of length $\rho$. In particular, they form up to reordering the reduced echelon basis of the vector subspace 
$$
V:=\left\langle v_1,\ldots,v_\rho\right\rangle\subset \kp^s
$$
that they generate over $\kp$ (in fact over any field). 
Hence, solving recombinations mainly reduces to find a system of $\kp$-linear equations that determine $V\subset \kp^s$. 

\medskip

Let $\mu=(\mu_1,\ldots,\mu_s)\in \kp^s$. Applying the logarithmic derivative with respect to $y$ to (\ref{prod}) and multiplying by $F$ we get 
\begin{equation}\label{maineq}
\mu\in V \iff \exists \,\alpha_1,\ldots,\alpha_s\in \kp\quad |\quad \sum_{i=1}^s \mu_i \hat{\F_i}\partial_y\F_i=\sum_{j=1}^s \alpha_j \hat{F_j}\partial_y F_j,
\end{equation}
with notations $\hat{F}_j=F/F_j$ and  $\hat{\F_i}=F/\F_i$.
The reverse implication holds since the $F_j$'s are supposed to be separable \cite[Lemma 1]{L1}. In \cite{L2}, the author show how to derive from (\ref{maineq}) a finite system of linear equations for $V$ that depends only on the $\F_i$'s truncated with $x$-adic precision $d_x+1$, \emph{assuming $F(0,y)$ separable of degree $d$}. For our purpose, we will rather consider $\vl$-adic truncation of the $\F_i$'s for a suitable $\lambda$, under the weaker hypothesis that $F$ is non degenerated.

\subsection{Residues and recombinations.}\label{ss2.2}
In what follows, we fix $\lambda=m/q\in \qp$. 
%
Given $G\in \kp((x))[y]$ and $\sigma\in \qp$, the $\vl$-truncation of $G$ with precision $\sigma$ is
$$
[G]_\lambda^{\sigma}:=\sum_{j+i\lambda\le \sigma} g_{ij} x^j y^i\in \kp[x^{\pm 1}][y].
$$
If $\lambda=0$, this is the classical Gauss (or $x$-adic) truncation $[G]_0^{\sigma}=G\mod x^{\sigma+1}$. If $G\in \kp[x,y]$, we can define the $\lambda$-degree of $G$,
$$
\dl(G):=\max(j+i\lambda,\, g_{ij}\ne 0).
$$
Note that $G=[G]_\lambda^{\dl(G)}$. Moreover, we have  $$\dl(GH)=\dl(G)+\dl(H)\quad {\rm and}\quad \dl(G+H)\le \dl(G)+\dl(H).$$
Let $\mu\in \fp^r$. Given the factorization \eqref{factoan}, we let
\begin{equation}\label{eq:Gmu}
G_{\mu}:=\sum_{i=1}^s \mu_i \big[\hat{\F_i}\partial_y\F_i\big]_\lambda^{\dl(F)}\,\, \in \,\,\kp[x,y].
\end{equation}
We denote by $\rho_k=\rho_k(\mu)$ the residues of $G_\mu/F$ at the roots $y_k\in \overline{\kp(x)}$ of $F$, that is
$$
\qquad \qquad \rho_{k}:=\frac{G_\mu(x,y_k)}{\partial_y F(x,y_k)}\,\, \in \,\,\overline{\kp(x)},\qquad k=1,\ldots,d=\deg_y(F).
$$
These residues are well defined since $F$ is separable. 
The next key result is mainly a consequence of \cite[Prop 8.7]{We16}:

\begin{prop}\label{prop:recombinations}
Suppose that $F$ is not degenerated. Then $\mu\in V$ if and only if $\rho_k\in \overline{\kp}$ for all $k=1,\ldots, d$.
\end{prop}

\begin{proof}
The direct implication follows from \eqref{maineq}. Let us prove the converse, assuming that the residues $\rho_k$ are constant. Let $\tau=\tau_{\lambda}$ as defined by \eqref{eq:tau_lambda}. Given $Q\in \kp((x))[y]$, we denote for short
$$
\tau_0(Q)=x^{-v_0(\tau(Q))} \tau(Q)\in \kp[[x]][y].
$$
Hence $\tau_0(F)\in \kp[x,y]$ is a primitive polynomial with primitive factors $\tau_0(F_j)$ in $\kp[x,y]$ and $\tau_0(\F_i)$ in $\kp[[x]][y]$.  Following \eqref{eq:v0vl}, we get
$$
n:=\deg_x(\tau_0(F))=q(\dl(F)-\vl(F)).
$$
The operators $\tau_0$ and $\partial_y$ commute and it's straightforward to check that
$$
\tau_0(G_\mu)=\sum_{i=1}^r \mu_i \big[\tau_0\widehat{(}\F_i)\partial_y\tau_0(\F_i)\big]_0^{n}.
$$
In other words, $\tau_0(G_\mu)$ coincides with the polynomial defined by  \eqref{eq:Gmu} when considering the recombinations of the analytic factors $\tau_0(\F_i)$ of $\tau_0(F)$ using the Gauss valuation $v_0$. 
Let $\phi_k(x):=x^{-m} y_k(x^q)$. We have
$
\tau_0(F)(x,\phi_k(x))=F(x^q,y_k(x^q))=0
$
for all $k$ so $\tau_0(F)$ has roots $\phi_1,\ldots,\phi_d$. Moreover,
$$\frac{\tau_0(G_\mu)(x,\phi_k(x))}{\partial_y \tau_0(F)(x,\phi_k(x))}=\frac{G_\mu(x^q,y_k(x^q))}{\partial_y F(x^q,y_k(x^q))}=\rho_{k}(x^q)\in \overline{\kp},
$$
so the residues of $\tau_0(G_\mu)/\tau_0(F)$ at the roots of $\tau_0(F)$ are constant by assumption. Since $F$ is separable and not degenerated, so is $\tau_0(F)$. Thus, we can apply \cite[Prop 8.7]{We16} to $\tau_0(F)$ and we deduce that $\tau_0(G_\mu)$ is a $\kp$-linear combination of the polynomials $\tau_0\widehat{(}F_j)\partial_y \tau_0(F_j)$, which in turns implies that $G_\mu$ is a $\kp$-linear combination of the polynomials $\hat{F_j}\partial_y F_j$. Hence $\mu \in V$ thanks to \eqref{maineq}, as required. 
\end{proof}

\begin{rem}
The assumption $F$ not degenerated is crucial to solve recombinations with $\vl$-precision $\dl(F)$. Otherwise, we might need to compute the $\F_i$'s with a higher precision. We refer the reader to \cite{We16} for various options to solve the recombination problem for arbitrary polynomials in the $x$-adic case.
\end{rem}

\subsection{Computing equations for $V$}\label{ss5.1}
Since $F$ is separable, $\rho_k$ belongs to the separable closure of $\kp(x)$ and we can talk about the derivative of $\rho_k$. Hence, an obvious necessary condition for that $\rho_k\in \overline{\kp}$ is that its derivative vanishes. More precisely, we have the following lemma:

\begin{lem}\label{constantresidue}
Let $p\ge 0$ be the characteristic of $\kp$. If $\rho_k'=0$ then $\rho_k\in \overline{\kp(x^p)}$. If moreover $p=0$ or $p\ge 2d (\dl(F)-\vl(F))$, then  $\rho_k\in \overline{\kp}$.
\end{lem}

\begin{proof}
If $\rho_k'=0$, then clearly $\rho_k\in \overline{\kp(x^p)}$. If $p=0$, the claim follows. If $p>0$, we consider the polynomial $\tau_0(F)$ defined above, of $x$-degree $n=q(\dl(F)-\vl(F))$. Its residue is $\rho_k(x^q)$ which thus lives in $\overline{\kp(x^{pq})}$. Hence, it's straightforward to check that we can divide the bound $p\ge 2 dn$  of \cite[Lemma 2.4]{GaoPDE} by $q$ in this context. 
\end{proof}

\noindent
Let us consider the $\kp$-linear operator 
\begin{equation}\label{eq:DGmu}
\begin{array}{ccc}
\D: \kp(x)[y] &\quad \longrightarrow &\quad \kp(x)[y] \\
G & \quad \longmapsto &  \quad \Big(G_x F_y-G_y F_x\Big)F_y-\Big(F_{xy}F_y-F_{yy} F_x\Big)G,
\end{array}
\end{equation}
with the standard notations $F_y$, $F_{xy}$, etc. for the partial derivatives. 

\begin{lem}\label{lem:rhoprimek}
We have $\rho_k'=0$ for all $k=1,\ldots,d$ if and only if $F$ divides $\D(G_\mu)$ in the ring $\kp(x)[y]$. 
\end{lem}

\begin{proof} 
Combining 
$
\rho_k(x)=\frac{G_\mu (x,y_k)}{F_y(x,y_k)}$ and $y_k'(x)=-\frac{F_x(x,y_k)}{F_y(x,y_k)},
$ 
we get 
$$
\rho_k'(x)=\frac{\D(G_\mu)(x,y_k)}{F_y^3(x,y_k)}.
$$
Thus $\rho_k'=0$ if and only if $\D(G_\mu)$ vanishes at all roots of $F$, seen as a polynomial in $y$.  The result follows since $F$ is separable.
\end{proof}

Let us denote $D_\mu:=\D(G_\mu)$ for short. We will need the following lemma:

\begin{lem}\label{lem:vlambdaDmu}
Suppose $\lambda\ge 0$. Then $3\vl(F)\le \vl(D_\mu)$ and $\dl(D_\mu)\le 3\dl(F)$.
\end{lem}

\begin{proof}
For any $Q\in \kp[x,y]$, the support of $x Q_x$ and $ yQ_y$ is contained in the support of $Q$. Hence $\vl(Q)\le \vl(x Q_x)$ and $\vl(Q)\le \vl(y Q_y)$ while $\dl(Q)\ge \dl(x Q_x)$ and $\dl(Q)\ge \dl(y Q_y)$. As $\vl(x)=\dl(x)=1$ and $\dl(y)=\vl(y)=\lambda\ge 0$, we get
\begin{equation}\label{eq:vlpartialy}
\quad \vl(Q_x),\,\,\vl(Q_y)\ge \vl(Q) \quad {\rm and}\quad  \dl(Q_x),\,\, \dl(Q_y)\le \dl(Q).
\end{equation}
In particular, we get from \eqref{eq:Gmu} that $\vl(G_\mu)\ge \vl(F)$. On the other hand we have $\dl(G_\mu)\le \dl(F)$ by the very definition \eqref{eq:Gmu}.  The claim then follows from \eqref{eq:DGmu}, using moreover that $\vl$ and $-\dl$ are valuations. 
\end{proof}

Lemma \ref{lem:rhoprimek} suggests to compute the $\vl$-adic euclidean division of $D_\mu$ by $F$ up to a sufficient precision to test divisibility in $\kp(x)[y]$. A difficulty is that $F$ is not necessarily $\lambda$-monic, hence we do not have access to Proposition \ref{prop:lambda_div}. To solve this issue, we adapt \cite[Section 5]{We16} to our context. We get: 

\begin{prop}\label{prop:final_division}
Given $\F_1,\ldots,\F_s$ with relative $\lambda$-precision $\dl(F)-\vl(F)$, we can compute a linear map 
$$
\phi:\kp^s\to \kp^{N}, \quad N \in \Oc( d(\dl(F)-\vl(F)))
$$ 
such that $\mu\in \ker(\phi)\iff F|D_\mu$, and so  with at most $\Ot(s N)$ operations in $\kp$. 
\end{prop}

\begin{proof}
Up to replace $F$ and $D(G_\mu)$ by their reciprocal polynomial, we may suppose that $\lambda\ge 0$ by Lemma \ref{lem:reciprocal}. Note first that $G_\mu$ only depends on the $\F_i$'s with relative $\lambda$-precision $\dl(F)-\vl(F)$ by \eqref{eq:vlpartialy} and Lemma \ref{lem:approx}. Let $0\le \alpha <q$ and $k\in \zp$ be the unique integers such that
$\tilde{F}:=\tau(x^k y^\alpha F)$ satisfies
\begin{equation}\label{eq:Ftilde}
q|\deg_y(\tilde{F})=d+\alpha \quad{\rm and}\quad 0\le v_0(\tilde{F})< q.
\end{equation} 
These conditions ensure that both $\tilde{F}$ and its leading coefficient $c:=\lc_y(\tilde{F})$ lie in the subring $\bpl\subset \kp[x,y]$ (Lemma \ref{lem:the_ring_apl}).  
Let $k'=k-2\vl(F)$ and $\tilde{D}_\mu:=\tau(x^{k'} y^{\alpha }D_\mu).$ 
Note that $F$ divides $D_\mu$ in $\kp(x)[y]$ if and only if $\tilde{F}$ divides $\tilde{D}_\mu$ in $\kp(x)[y]$. Moreover, $k'$ does not depend on $\mu$ so the map $\mu\mapsto \tilde{D}_\mu$ is $\kp$-linear. 

\medskip
\noindent
\textbf{Claim.} We have $v_0(\tilde{D}_\mu)\ge v_0(\tilde{F})$ and $\deg_x(\tilde{D}_\mu)\le 3\deg_x(\tilde{F})$. 

\medskip
\noindent 
\emph{Proof of the claim.} As $\lambda\ge 0$, Lemma \ref{lem:vlambdaDmu} and $v_0(\tau(Q))=q \vl(Q)$  give
$$
v_0(\tilde{D}_\mu)=q(k'+\alpha\lambda +\vl(D_\mu)\ge q(k+\alpha\lambda +\vl(F))=v_0(\tilde{F})\ge 0.
$$
In a similar way, using now $\deg_x(\tau(Q))=q\dl(Q)$, we get
\begin{eqnarray*}
\deg_x(\tilde{D}_\mu) &\le & q(k'+\alpha\lambda + 3\dl(F)) \\ &=& 2q(\dl(F)-\vl(F))+q(k+\alpha\lambda +\vl(F)) \\
 &=& 2(\deg_x(\tilde{F})-v_0(\tilde{F}))+\deg_x(\tilde{F})  \\
 &\le & 3\deg_x(\tilde{F}).  \qquad\qquad\qquad\qquad\hfill\square
\end{eqnarray*}
\noindent 
As $0\le v_0(\tilde{F}) \le v_0(\tilde{D}_\mu)$, $\tilde{F}$ divides $\tilde{D}_\mu$ in $\kp(x)[y]$ if and only if it divides $\tilde{D}_\mu$ in $\kp[x][y]$ by Gauss Lemma. To reduce to the monic case, we localize $\kp[x]$ at some prime $a\in \kp[x^q]$ coprime to $c:=\lc_y(\tilde{F})$. The euclidean division 
\begin{equation}\label{eq:euclidean_division}
\tilde{D}_\mu=Q_{\mu}\tilde{F}+R_\mu \in \kp[x]_{(a)}[y]
\end{equation}
is now well defined. Any $Q\in \bpl\subset \kp[x,y]$ has a unique $a$-adic expansion
\begin{equation}\label{eq:a-adic}
Q=\sum_{i=0}^{\lfloor \deg(Q)/\deg(a)\rfloor} q_i(x,y) a(x)^i,\quad {\rm with}\quad  q_i\in \bpl \quad {\rm and}\quad \deg_x q_i < \deg\,a.
\end{equation}
Note that $q_i\in \bpl$ since $a\in \bpl$. Let $\big\{Q\big\}^{n}_m=\sum_{i=m}^{n} q_i a^i$ and $\big\{Q\big\}^{n}=\big\{Q\big\}^{n}_0$.
Since $\deg_x(\tilde{D}_\mu)\le 3 \deg_x(\tilde{F})$, we deduce from (the proof of) \cite[Lemma 5.2]{We16}  that $\tilde{F}$ divides $\tilde{G}$ if and only if 
$$\big\{Q_{\mu}\big\}^{n}_{m}=\big\{R_\mu \big\}^{n}=0,\quad {\rm with}\quad 
m:=\Big\lfloor \frac{2d_x}{\deg\,a} \Big\rfloor+1 \quad{\rm and} \quad n:=\Big\lceil \frac{3d_x}{\deg\,a}\Big\rceil,
$$
where $d_x=\deg_x(\tilde{F})$. We have $d_x=q(\dl(F)-\vl(F))$ by $\eqref{eq:Ftilde}$. Since both polynomials $\{Q_\mu\}^n_m$ and $\{R_\mu\}^n$ live in $\bpl$, we deduce that their supports have size $\Oc(d(\dl(F)-\vl(F)))$. The linear map $$\phi(\mu):=\left(\{Q_\mu\}^n_m/a^m, \{R_\mu\}^n\right)$$ thus satisfies the conditions of Proposition \ref{prop:final_division}.  Let us look at complexity issues. If $Q_1,Q_2\in \bpl$ have $x$-degrees $\Oc(d_x)$ and relative $y$-degrees $\deg_y(Q_i)-v_y(Q_i)\in \Oc(d)$, we compute $\{Q_1\}^n$, $\{Q_2\}^n$ and $\{Q_1 Q_2\}^n$ in time $\Ot(d d_x/q)$ thanks to Proposition \ref{prop:mult} since all operations in \eqref{eq:a-adic} take place in $\bpl$. We have $c\in \kp[x^q]\subset \bpl$ invertible modulo $a$, and computing $\{c^{-1}\}^n$ costs $\Ot(d_x/q)$. Then, adapting  the proof of Proposition \ref{prop:div} in the $a$-adic case, we compute \eqref{eq:euclidean_division} with $a$-adic precision $n$ and thus $\phi(\mu)$ in time $\Ot( d d_x/q)=\Ot(d(\dl(F)-\vl(F)))$. To build the matrix of $\phi$, we compute $\phi(\mu_i)$ where the $\mu_i$'s run over the canonical basis of $\kp^s$.  Given the $\F_i$'s with relative $\lambda$-precision $\dl(F)-\vl(F)$, computing $G_{\mu_i}=\big[\hat{\F_i}\partial_y\F_i\big]_\lambda^{\dl(F)}$  costs $\Ot(d(\dl(F)-\vl(F))$ thanks to Corollary \ref{cor:lambda_mult}. Summing over all $i=1,\ldots,s$, we get the result. 
\end{proof}

\begin{rem}
We need to compute $a\in \kp[x^q]$ coprime to $c$. As $a=a_0(x^q)$ and $c=c_0(x^q)$, we look for $a_0$ coprime to $c_0$. We have $\deg_x(c_0)\le d_x/q=\dl(F)-\vl(F)$. If $\Card(\kp)\ge \dl(F)-\vl(F)$, we use multipoint evaluation of $c_0$ at $\deg_x(c_0)$ distinct elements of $\kp$ to find $z\in \kp$ such that $c_0(z)\ne 0$, and we take $a(x)=x^q-z$. Otherwise, we follow a similar strategy in a finite extension of $\kp$, considering now $a=a_0(x^q)$, with $a_0$ the minimal polynomial of $z$ over $\kp$. The cost fits in the aimed bound.
\end{rem}

\begin{cor}\label{cor2}
If $\kp$ has characteristic zero or greater than $2d(\dl(F)-\vl(F))$ and $F$ is non degenerated, then $(v_1,\ldots, v_\rho)$ is the reduced echelon basis of $\ker(\phi)$. 
\end{cor}

\begin{proof}
Follows from Proposition \ref{prop:recombinations}, Lemma \ref{constantresidue}, Lemma \ref{lem:rhoprimek} and Proposition \ref{prop:final_division}.
\end{proof}

If $\kp$ has small characteristic $p$, we need extra conditions to ensure $\rho_k\in \bar{\kp}$. These conditions rely on linear algebra over the prime field $\fp_p$ of $\kp$. They are based on Niederreiter's operator, which was originally introduced for univariate factorization over finite fields \cite{Nied}, and used then for bivariate factorization in \cite{L2}. We deliberately do not go into the details here. We assume $\lambda\ge 0$. We introduce the following $\fp_p$-linear map:
\begin{eqnarray*}
\psi:\ker(\phi_{|\fp_p^s}) &\longrightarrow &\kp[x^p,y^p]_{p\dl,p(d-1)}\\
\mu &\longmapsto & G_{\mu}^p-\partial_y^{p-1}(G_\mu F^{p-1}).
\end{eqnarray*}
In contrast to \cite{L2}, the subscripts indicate the $\lambda$-degree and the $y$-degree.

\begin{prop}\label{prop:small_caracteristic}
The map $\psi$ is well-defined and $(v_1,\ldots, v_\rho)$ is the reduced echelon basis of $\ker(\psi)$. 
\end{prop}

\begin{proof}
We check that $\partial_y(\psi(\mu))=0$, so $\psi(\mu)$ is a polynomial in $y^p$ of $y$-degree $p(d-1)$. Since $\dl(Q_y)\le \dl(Q)$ (proof of Lemma \ref{lem:vlambdaDmu}), $\psi(\mu)$ has $\lambda$-degree at most $p\dl$. Since moreover $\mu\in \ker(\phi)$, we have $\rho_k\in \overline{\kp(x^p)}$ by Lemma \ref{constantresidue} and Proposition \ref{prop:final_division}, which forces $\psi(\mu)$ to be a polynomial in $x^p$ (see \cite[Lemma 4]{L2}). Hence $\psi$ is well-defined. If $\lambda=0$, the second claim follows from \cite[Proposition 2]{L2} . If $\lambda\ne 0$,  
we reason as in the proof of Proposition \ref{prop:recombinations}, passing through the polynomials $\tau_0(F)$ and $\tau_0(G_\mu)$ to reduce to the case $\lambda=0$ (using again that $\partial_y$ and $\tau_0$ commute). 
\end{proof}

\begin{prop}\label{prop:solving_recombination}
Denote $N=d (\dl(F)-\vl(F))$. Assume $F\in \kp[x,y]$ non degenerated. Given $\F_1,\ldots,\F_s$ with relative $\lambda$-precision $\dl(F)-\vl(F)$, we can solve the recombination problem with
\begin{enumerate}
\item $\Ot(sN)+\Oc(s^{\omega-1}N)$ operations in $\kp$ if $p=0$ or $p\ge 2N$,
\item $\Oc(k s^{\omega-1}N)$ operations in $\fp_p$ if $\kp=\fp_{p^k}$. 
\end{enumerate}
\end{prop}

\begin{proof}
We can compute the reduced echelon basis of the kernel of a matrix of size $s\times N$ with coefficient in a field $\lp$ with $\Oc(s^{\omega-1} N)$ operations in $\lp$ \cite[Theorem 2.10]{Stor}. Hence, the first point follows from Proposition \ref{prop:recombinations} and Corollary \ref{cor2}. Suppose that $\kp=\fp_{p^k}$. Thus $\kp$ is an $\fp_p$-vector space of dimension $k$ and it follows again from Proposition \ref{prop:recombinations} that we can build the matrix of $\phi_{|\fp_p^s}$ and compute a basis of its kernel over $\fp_p$ in the aimed cost. To build the matrix of $\psi$ we reason again with the polynomials $\tau_0(F)$ and $\tau_0(G_\mu)$ to reduce to the case $\lambda=0$, using that $\partial_y$ and $\tau_0$ commute. We apply then \cite[Proposition 13]{L2}, using again that the complexity can be divided by $q$ since we work in the sparse subring $\bpl\subset \kp[x,y]$ (in the non monic case, we localize at some $a\in \kp[x^q]$ as in the proof of Proposition \ref{prop:recombinations}). The resulting complexity fits in the aimed cost. The matrix of $\psi$ having size at most $s\times k N$ over $\fp_p$, we conclude.
\end{proof}

\subsection{Proof of Theorem \ref{thm:main} and Corollary \ref{cor:mmain}} \label{smallchar}
The key point is to choose a good slope $\lambda$ before applying Proposition \ref{prop:solving_recombination}. Let $F\in \kp[x,y]$ of $y$-degree $d$ with Newton polygon $N(F)$ and lower convex hull $\Lambda(F)$. Let $V=\Vol(N(F))$.

\begin{lem}\label{lem:good slope}
Let $\lambda:=\lambda_F$ be the average slope of $\Lambda(F)$ (Definition \ref{def:average_slope}). Assume that $y$ does not divide $F$. Then
$$V\le d(\dl(F)-\vl(F))\le 2 V.$$
\end{lem}

\begin{proof}
It is a similar proof as that of Proposition \ref{prop:key_mlG}. Consider the bounding parallelogram $ABCD$ of $N(F)$ with two vertical sides and two sides of slope $-\lambda$ (figure \ref{fig good slope} below). We have 
$\Vol(ABCD)=d(\dl(F)-\vl(F))$ which gives the first inequality.
Consider $I$ and $J$ the left and right end points of $\Lambda(F)$ and let $K\in [BC]\cap N(F)$ and $L\in [AD]\cap N(F)$. Then 
$$
V\ge \Vol(IJK)+\Vol(IJL)=\frac{\Vol(IBCJ)}2 + \frac{\Vol(IADJ)}2=\frac{\Vol(ABCD)}2,
$$
the inequality since $IJK$ and $IJL$ are contained in $N(F)$, and the first equality since $(IJ)$ is parallel to $(AD)$ and $(CD)$ by choice of $\lambda=\lambda_F$. The result follows.
\end{proof}

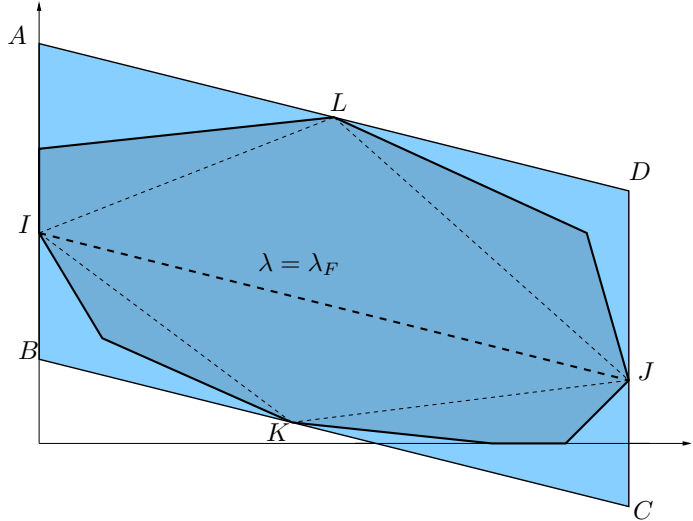
\begin{figure}[h]
	\caption{Proof of Lemma \ref{lem:good slope}. In dark blue the polygon $N(F)$ and in light blue its bounding parallelogram of slope $\lambda_F$.}
	\label{fig good slope}
	\bigskip
	\scalebox{0.28}{\begin{picture}(0,0)%
\includegraphics{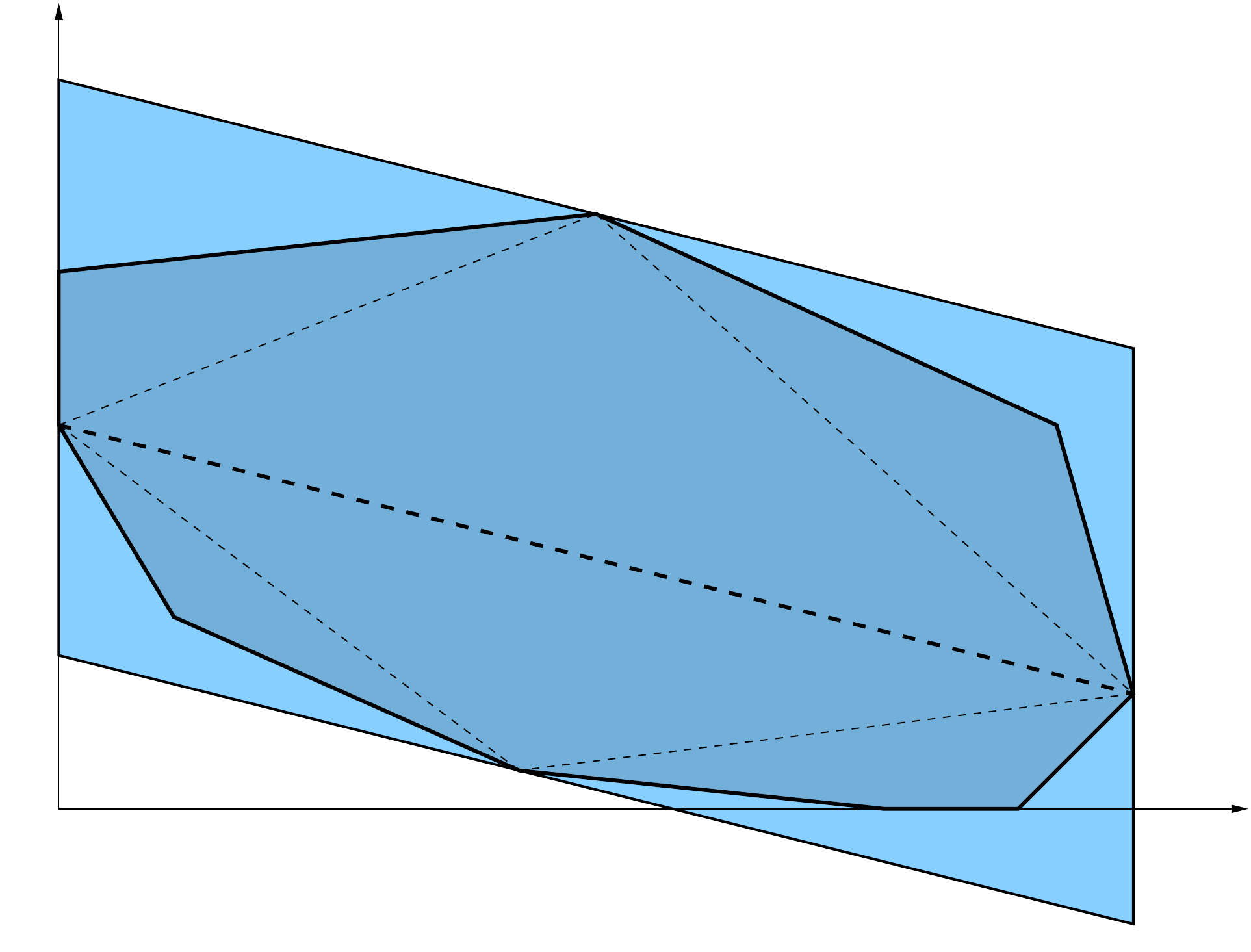}%
\end{picture}%
\setlength{\unitlength}{4144sp}%
\begingroup\makeatletter\ifx\SetFigFont\undefined%
\gdef\SetFigFont#1#2#3#4#5{%
  \reset@font\fontsize{#1}{#2pt}%
  \fontfamily{#3}\fontseries{#4}\fontshape{#5}%
  \selectfont}%
\fi\endgroup%
\begin{picture}(14662,11147)(2461,-12536)
\put(2476,-2311){\scalebox{3.5}{$A$}}
\put(2701,-9061){\scalebox{3.5}{$B$}}
\put(15751,-5236){\scalebox{3.5}{$D$}}
\put(2701,-6361){\scalebox{3.5}{$I$}}
\put(9361,-3751){\scalebox{3.5}{$L$}}
\put(7831,-7171){\scalebox{3.5}{$\lambda=\lambda_F$}}
\put(8000,-10800){\scalebox{3.5}{$K$}}
\put(15931,-9511){\scalebox{3.5}{$J$}}
\put(15841,-12481){\scalebox{3.5}{$C$}}
\end{picture}%
}
\end{figure}

\medskip
\noindent
Previous results lead to algorithm \texttt{Factorization} below.

\begin{algorithm}[ht]
  \nonl\TitleOfAlgo{\texttt{Factorization}($F$)}
  \KwIn{%
    $F\in \kp[x,y]$ primitive, separable in $y$ and non degenerated.}%
  \KwOut{The irreducible factorization of $F$ over $\kp$}%
  \leIf{$y$ divides $F$}{$L=[y]$ and $F\gets F/y$}{$L\gets [\,\,]$}%
 $\lambda\gets \lambda_F$ and $\sigma\gets \dl(F)-\vl(F)+\ml(F)$\;%
 $[\F_1,\ldots,\F_s] \gets$ \texttt{Facto}($F,\lambda,\sigma$)\;%
  \lIf{$s=1$}{\Return{$L$}}%
 Compute the reduced echelon basis $(v_1,\ldots,v_\rho)$ of $V$ using Proposition \ref{prop:solving_recombination}\;%
\For{$j=1,\ldots,\rho$}{
Compute $\tilde{F}_j:=[\lc_y(F)\prod_{i=1}^s \F_i^{v_{ji}}]_\lambda^{\dl(F)}$\;%
Compute the primitive part $F_j$ of $\tilde{F}_j$ with respect to $y$}%
\Return{$L\cup [F_1,\ldots, F_\rho]$}
\end{algorithm}

\begin{prop}\label{prop:algo_factors}
Algorithm \emph{\texttt{Factorization}} is correct.
Up to the cost of univariate factorizations, it takes at most
\begin{enumerate}
\item $\Ot(s V)+\Oc(s^{\omega-1}V)$ operations in $\kp$ if $p=0$ or $p\ge 4V$,
\item $\Oc(k s^{\omega-1}V)$ operations in $\fp_p$ if $\kp=\fp_{p^k}$.
\end{enumerate}
\end{prop}

\begin{proof}
By Theorem \ref{prop:fast_facto}, Step 3 computes the $\F_i$'s with relative $\lambda$-precision at least $\dl(F)-\vl(F)$. Thus, Proposition \ref{prop:solving_recombination} and Lemma \ref{lem:good slope} ensure that the $v_j$'s at step 5 are solutions to the recombination problem \eqref{eq:eq2}. Since $F$ is primitive, so are the $F_j$'s. Since $\dl(\lc(F)/\lc(F_j)) + \dl(F_j ) \le \dl(F)$ we have $\tilde{F}_j= \frac{\lc(F)}{\lc(F_j)}F_j$ so $F_j$ is the primitive part of $\tilde{F}_j$. Hence the algorithm returns a correct answer. Since $\ml(F)\le \dl(F)-\vl(F)$ (Definition \ref{def:mlambda}), we have $\sigma\le 4V/d$ by Lemma \ref{lem:good slope}. Hence step 3 costs $\Ot(V)$ by Theorem \ref{prop:fast_facto}. Step 5 fits in the aimed bound by Proposition \ref{prop:solving_recombination} and Lemma \ref{lem:good slope}. Using technique of subproduct trees, step 7 costs $\Ot(d(\dl(F)-\vl(F))=\Ot(V)$ by Corollary \ref{cor:lambda_mult}, and computing primitive parts at step 8 has the same cost. This concludes the proof.
\end{proof}

\subsubsection*{Proof of Theorem \ref{thm:main}.}
It follows immediately from Proposition \ref{prop:algo_factors} since $s$ is smaller or equal to the lower lattice length $r$ of $N(F)$. Note that testing non degeneracy amounts to test squarefreeness of some univariate polynomials whose degree sum is $r$, hence costs only $\Ot(r)$ operations in $\kp$. $\hfill\square$

\subsubsection*{Proof of Corollary \ref{cor:mmain}.}
The corollary follows straightforwardly from Theorem \ref{thm:main}. However, let us explain for the sake of completeness how to compute quickly the minimal lower lattice length $r_{0}(F)=r_{0}(N(F))$. Recall from  \eqref{eq:rmin} that for a lattice polygon $P$, 
$$
r_{0}(P)=\min\{r(\tau(P)),\, \tau \in \Aut(\zp^2)\}
$$
where $r(\tau(P))$ stands for the lattice length of the lower convex hull $\Lambda(\tau(P))$, and $\Aut(\zp^2)$ stands for the group of affine automorhisms. 

\begin{lem}\label{lem:rmin}
Let $P$ be a lattice polygon, with edges $E_1,\ldots, E_{n}$. Denote $w_i\in \zp^2$ the inward orthogonal primitive vector of $E_i$. There exist $\tau_i,\tau_i' \in GL_2(\zp)$ with $\det(\tau_i)=1$ and $\det(\tau_i')=-1$ and such that $\tau_i(w_i)=\tau_i'(w_i)=(1,0)$. Then
$$
r_{0}(P)=\min\left(r(\tau_1(P)),r(\tau_1'(P)), \ldots, r(\tau_n(P)),r(\tau_n'(P))\right). 
$$
Geometrically, the maps $\tau_i$ and $\tau'_i$ simply send $E_i$ to a vertical left hand edge. Such maps are straightforward to compute (note that they are not unique). 
\end{lem}

\begin{proof}
Since the lower lattice length is invariant by translation, it's sufficient to look for a map $\tau\in GL_2(\zp)$ that reaches $r_0$. Let us first consider $\tau\in GL_2(\rp)$. Consider the set $I_\tau=\{j,\, \tau(E_j)\subset \Lambda(\tau(P)\}$ of the indices of the lower edges of $\tau(P)$. Denoting $d_j(\tau)=\det((1,0),\tau(w_j))$, we have
$$
j\in I_\tau \iff d_j(\tau)>0.
$$
The maps $\tau\mapsto d_j(\tau)$ being continuous, we deduce that $I_\tau\subset I_{\tau'}$ for all $\tau'\in GL_2(\rp)$ close enough to $\tau$, and with equality $I_\tau= I_{\tau'}$ if $d_j(\tau)\ne 0$ for all $j=1,\ldots,n$. Obviously, if $\tau,\tau'\in GL_2(\zp)$ then $I_\tau\subset I_{\tau'}$ implies $r(\tau(P))\le r(\tau'(P))$ and equality $I_\tau= I_{\tau'}$ implies equality of the lower lattice lengths. It follows that $r_0$ is reached at $\tau\in GL_2(\zp)$ such that $d_i(\tau)=0$ for some $i$ (such a $\tau$ exists for each $i$). This forces $\tau(w_i)=\pm (1,0)$ and we may suppose $\tau(w_i)=(1,0)$ since the lower lattice length is invariant by vertical axis symmetry. But if $\tau'\in GL_2(\zp)$ is another map such that $\tau'(w_i)=(1,0)$ and which satisfies moreover $\det(\tau)=\det(\tau')$, then 
$$
d_j(\tau')=\det(\tau'(w_i),\tau'(w_j))=\det(\tau')\det(w_i,w_j)=\det(\tau)\det(w_i,w_j)=d_j(\tau)
$$
for all $j=1,\ldots,n$, from which it follows that $I_\tau=I_{\tau'}$, hence $r(\tau(P))=r(\tau'(P))$. The lemma follows.
\end{proof}

\bibliographystyle{abbrv}{\bibliography{tout}}

\end{document}